\def\versiondate{18 June 2003}
\input math.macros
\input Ref.macros

%\proofmodetrue
%\leftsectionheadtrue
\checkdefinedreferencetrue
%\continuousnumberingtrue
\continuousfigurenumberingtrue
\theoremcountingtrue
\sectionnumberstrue
%\figuresectionnumberstrue
\forwardreferencetrue
%\lefteqnumberstrue
\citationgenerationtrue
\nobracketcittrue
\tocgenerationtrue
\hyperstrue
\initialeqmacro

\input bases.key

\bibsty{myapalike}

\def\beginbullets{\begingroup\everypar={$\bullet$\enspace}}
\def\endbullets{\smallbreak\endgroup}
\def\fusf{{\ss FSF}}
\def\fsf{{\ss FSF}}
\def\wusf{{\ss WSF}}
\def\wsf{{\ss WSF}}

\def\B{{\cal B}}
\font\frak=eufm10   %% or  eufb10
\font\scriptfrak=eufm7
\font\scriptscriptfrak=eufm5
\def\mathfrak#1#2{%       %% This cannot be done as
\def#1{{\mathchoice%
{{\hbox{\frak #2}}}%
{{\hbox{\frak #2}}}%
{{\hbox{\scriptfrak #2}}}%
{{\hbox{\scriptscriptfrak #2}}}}}}
\mathfrak{\ba}{B}  %% random variable for base
\mathfrak{\qba}{S}  %% random variable for subset
\mathfrak{\fo}{F}  %% random forest
\def\M{{\cal M}}

\def\F{{\cal F}}
\def\Fd{{\Bbb F}}  %% field

\def\sgn{\mathop{\rm sgn}}

\def\C{{\Bbb C}}
\def\H{{\ell^2(E)}}
\def\bd{{\partial}}
\def\Ext{{\rm Ext}}
\def\u{{\bf u}}
\def\v{{\bf v}}
\def\w{{\bf w}}
\def\STAR{\bigstar}        %% The space spanned by the stars.
\def\CYCLE{\diamondsuit}
          %% A star.
\def\vertex{{\ss V}}
\def\edge{{\ss E}}
\def\edges{{\ss E}}
\def\setminu{\setminus}
\def\id{{\rm id}}  %% The identity operator.
\def\ev#1{{\cal #1}}
\def\A{{\ev A}}
\def\ip#1{(\changecomma #1)}
\def\bigip#1{\big(\bigchangecomma #1\big)}
\def\biggip#1{\biggl(\bigchangecomma #1\biggr)}
\def\leftip#1{\left(\leftchangecomma #1\right)}  
\def\changecomma#1,{#1,\,}
\def\bigchangecomma#1,{#1,\;}
\def\leftchangecomma#1,{#1,\ }
\def\d#1{#1^\dagger}

\def\CHI{\mathchoice%
{\raise2pt\hbox{$\chi$}}%
{\raise2pt\hbox{$\chi$}}%
{\raise1.3pt\hbox{$\scriptstyle\chi$}}%
{\raise0.8pt\hbox{$\scriptscriptstyle\chi$}}}
\def\Var{{\rm Var}}
\def\CC{{M}}
\def\HFS{{H^F_S}}
\def\HSFS{{H_{S, F \setminus S}}}
\def\D{{\cal D}}  %% dual of multivector (*?)
\def\QF{{\widetilde Q_{F, B}}}
\def\pd#1{\P^\bd_{#1}}  %% probability measure on #1-subcomplexes 
\def\pcx#1{\P^{\rm cx}_{#1}}  %% probability measure on independent #1-sets
\def\SOTto{\buildrel {\rm SOT} \over \longrightarrow}
\def\T{{\Bbb T}}  % circle
\def\ent{{\ss Ent}}  %% entropy
\def\constant#1{{\bf #1}}
\def\asp{\alpha}  % aspect
\def\bm{{\ss BM}}  % B-M density
%for Okounkov's example:
\def\la{\lambda}
\mathfrak{\fM}{M}
\mathfrak{\fS}{S}
\def\KK{{\ss K}}
\def\frac#1#2{{#1 \over #2}}
%
%for BOO:
\def\Fr{{\rm Fr}}
\def\ix{\theta}
\def\KK{{\ss K}}
\def\kk{{\ss k}}
\def\sh{{1 \over 2}}

\def\BLPSgip{\ref b.BLPSgip/}%\def\BLPSgip{\htmllocref{BLPSgip}{[BLPS99]}}}
\def\BLPSusf{\ref b.BLPSusf/, hereinafter referred to as BLPS (2001)%
\def\BLPSusf{BLPS \htmllocref{BLPSusf}{(2001)}}}

\def\firstheader{\eightpoint\ss\underbar{\raise2pt\line 
    {
    %{To appear in {\it Ann.\ Probab.}
    \hfil Version of \versiondate}}}

\beginniceheadline

%\ifproofmode \relax \else\head{} {Version of \versiondate}\fi 
\vglue20pt

\title{Determinantal Probability Measures}

%\bigskip
%\centerline{\tt NOT FOR DISTRIBUTION}
%\bigskip

\author{Russell Lyons}

\abstract{%
Determinantal point processes have arisen in diverse settings in recent
years and have been investigated intensively.
We study basic combinatorial and probabilistic aspects
in the discrete case.
Our main results concern relationships with matroids, stochastic
domination, negative association, completeness for infinite matroids, tail
triviality, and a method for extension of results from orthogonal
projections to positive contractions.
We also present several new 
avenues for further investigation, involving
Hilbert spaces, combinatorics, homology, and group representations,
among other areas.
}

\bottomII{
Primary 60B99, % probability thy on algebraic and topological structures
05B35.  %  Matroids
Secondary
60D05, % geometric probability, random sets
05C05, %Trees
15A75.} % Exterior algebra, Grassmann algebras.
{Spanning trees, matroids, exterior algebra, probability, determinants,
orthogonal projections, positive contractions, negative association, bases.}
{Research partially supported by a Varon Visiting Professorship at
the Weizmann Institute of Science, NSF grants DMS-9802663 and
DMS-0103897, Microsoft Corp., and the Miller Institute for Basic Research in
Science at the University of Calif., Berkeley.}

%\vfill\eject
\articletoc

\bsection {Introduction}{s.intro}

A determinantal probability measure is one whose elementary
cylinder probabilities are given by determinants.
More specifically, suppose that $E$ is a finite or countable set and that
$Q$ is an $E \times E$ matrix.
For a subset $A \subseteq E$, let $Q\restrict A$ denote the submatrix of
$Q$ whose rows and columns are indexed by $A$.
If $\qba$ is a random subset of $E$ with the property that for all finite $A
\subseteq E$, we have 
$$
\P[A \subseteq \qba] = \det (Q\restrict A)
\,,
\label e.DPM
$$
then we call $\P$ a {\bf determinantal probability measure}.
For a trivial example, if $\P$ is the Bernoulli($p$) process on $E$, then \ref
e.DPM/ holds with $Q = pI$, where $0 \le p \le 1$ and $I$ is the identity
matrix.
Slightly more generally, when $\P$ is a product measure, there is a diagonal
matrix $Q$ that makes \ref e.DPM/ true. 

The inclusion-exclusion principle in combination with \ref e.DPM/
determines the probability of each elementary
cylinder event (as a determinant, in
fact; see \ref r.incexc/). Therefore,
for every $Q$, there is at most one probability
measure satisfying \ref e.DPM/.
Conversely, it is known (see, e.g., \ref s.contract/) that there is a
determinantal probability measure corresponding to $Q$ if $Q$ is the matrix of
a positive contraction on $\ell^2(E)$.

The continuous analogue of determinantal probability measures,
i.e., determinantal point processes in $\R^n$, have a long history.
This began in the 1950s with Wigner's investigation of the distribution of
eigenvalues of certain large random matrices in order to study energy levels of
large atoms.
The study of eigenvalues of random matrices continues to be an important
topic in physics.
\ref b.Dyson/ proved that the so-called correlation functions
of the distributions that Wigner considered could be described by simple
determinants.
In the early 1970s, a conversation with Dyson led
Montgomery to realize that conjectures he was then formulating concerning
the zeroes of the Riemann zeta function were related to the distribution of
eigenvalues of random matrices.
This idea has been extremely fruitful; see, e.g., the recent reviews by \ref
b.Conrey/ and \ref b.Diaconis:Gibbs/.
At about that same time, in studying fermions in physics, 
general point processes with a determinantal form were introduced by
\ref b.Macchi/ (see also the references therein) and are known
as {\bf fermionic processes} in mathematical physics.
The discrete case, which is the one studied here,
first appeared in
Exercises 5.4.7--5.4.8 of the book by \ref b.DVJ/, where it is noted that the
continuous point processes can be obtained as scaling limits of the discrete
measures.
Since the end of the 20th century, both the continuous and discrete
measures have received much attention, especially, specific measures and, if
discrete, their continuous scaling limits.
For example, scaling limits of discrete instances
have been studied in \ref b.Bor:char/,
\refbmulti{BO:hyperg,BO:z-measures,BO:harmonic},
\ref b.BOO:plancherel/,
\refbmulti{Joh:DisOP,Joh:noninter}, \ref b.Oko:wedge/,
and \ref b.OkoResh:schur/.
The literature on the purely continuous case is too voluminous to list here.
A general study of the discrete and continuous cases has been
undertaken independently by
\refbmulti{ShiTak:I,ShiTak:II} (announced in \ref b.ShiTak:ann/) and
\ref b.ShiYoo:glauber/, but there is little overlap with our work here.
Several aspects of general stationary determinantal probability measures on
$\Z^d$ are studied by \ref b.LS:dyn/.
For a survey of both the discrete and continuous cases, see
\ref b.Soshnikov:survey/.

Our purpose is to establish some new basic combinatorial and probabilistic
properties of all (discrete) determinantal probability measures.
However, for the benefit of the reader who has not seen any before, 
we first display a few examples of such measures.
Most examples in the literature are too involved even to state here.
We restrict our examples to a few that can be detailed easily.

\procl x.ust
The most well-known example of a (nontrivial discrete) determinantal
probability measure is that where $\qba$ is a uniformly chosen random spanning
tree of a finite connected graph $G = (\vertex, \edge)$ with $E := \edge$.
In this case, $Q$ is the {\bf transfer current matrix} $Y$,
which is defined as follows.
Orient the edges of $G$ arbitrarily.
Regard $G$ as an electrical network with each edge having unit
conductance.  Then $Y(e, f)$ is the amount of current flowing along the edge
$f$ when a battery is hooked up between the endpoints of $e$ of such
voltage that in the network as a whole, unit current flows from the tail of
$e$ to the head of $e$.
The fact that \ref e.DPM/ holds for the uniform spanning tree is due to
\ref b.BurPem/ and is called the Transfer Current Theorem.
The case with $|A| = 1$ was shown much earlier by \ref b.Kirchhoff/, while
the case with $|A| = 2$ was first shown by \ref b.BSST/.
\endprocl

\procl x.renewal Let $0 < a < 1$ and define the Toeplitz matrix
$$
R(i, j) := {1-a \over 1+a} a^{|j - i|}
$$
for $i, j \in \Z$.
Then $\P^R$ is a renewal process on $\Z$ (\ref B.Soshnikov:survey/) with
renewals at each point of $\qba$.
The distance between successive points of $\qba$ has the same distribution as
one more than the number of tails until 2 heads appear for a coin that has
probability $a$ of coming up tails.
\ref b.LS:dyn/ extend this example to other regenerative processes.
\endprocl

\procl x.BOO
For $\ix>0$, consider 
the probability measure on the set of all partitions $\la$ of all nonnegative
integers
$$
M^\ix (\la) :=
e^{-\ix} \ix^{|\la|} \,
\left(\frac{\dim\la}{|\la|!}
\right)^2 \,. 
$$
Here, $|\la|$ is the sum of the entries of $\la$ and $\dim\la$ is the
dimension of the representation of the symmetric group on $|\la|$ letters that
is determined by $\la$.
This measure is a Poissonized version of the Plancherel measure.
In order to exhibit $M^\ix$ as a determinantal probability measure, as
established
by \ref b.BOO:plancherel/, we first represent partitions as subsets of
$\Z+\sh$ as follows.
Consider a partition $\la$ as a Young diagram.
Call $d$ the number of squares on the diagonal of $\la$.
The Frobenius
coordinates of $\la$ are $ \{ p_1, \ldots, p_d, q_1, \ldots, q_d \}$,
where $p_i$ is the number of
squares in the $i$th row  to the right of the diagonal
and $q_i$ is number of squares in the $i$th column
below the diagonal. 
Following \ref b.VK/, define  the {\bf modified 
Frobenius coordinates} $\Fr(\la)$ of $\la$ by 
$$
\Fr(\la):=
\left\{p_1+\sh,\dots,p_d+\sh,-q_1-\sh,\dots,-q_d-\sh
\right\} \subset \Z+\sh\,.
$$
Clearly, the map $\la \mapsto \Fr(\la)$ is injective.
Then the law of $\Fr(\la)$ is $\P^\KK$ when $\la$ has the law $M^\ix$,
where the matrix $\KK$ is 
$$
\KK(x,y):=\cases{
\displaystyle
\sqrt{\ix}\,\,\frac{\kk_+(|x|,|y|)}{|x|-|y|} &if $xy >0$,\cr
\displaystyle \sqrt{\ix}\,\,\frac{\kk_-(|x|,|y|)}{x-y} &if $xy<0$,\cr
}$$
with the functions $\kk_\pm$ defined by
\begineqalno
\kk_+(x,y)&:= 
J_{x-\sh}\, J_{y+\sh}  - J_{x+\sh} \, J_{y-\sh} \,, \cr
\kk_-(x,y)&:= J_{x-\sh}\, J_{y-\sh}  + J_{x+\sh} \, J_{y+\sh} \,,
\endeqalno
where $J_x:=J_x(2\sqrt{\ix})$ is the Bessel function of
order $x$ and argument $2\sqrt{\ix}$. 
Although $\KK$ is not self-adjoint on $\ell^2(\Z+\sh)$, its restrictions to
$\ell^2(\Z^+-\sh)$ and to $\ell^2(\Z^-+\sh)$ are.
\endprocl

\comment{
Okounkov Schur measures: 

Consider the set of all partitions $\la$ and introduce
the following
function of $\la$ 
$$
\fM(\la)={1 \over Z} |s_\la(x)|^2
\,,
$$
where $s_\la$ are the Schur functions (\ref b.Macd/)
in auxiliary variables $x_1,x_2,\dots$
$$
Z:= \sum_\la |s_\la(x)|^2 = \prod_{i,j} (1-x_i \overline{x_j})^{-1} \,.
$$
$\fM$ is a probability measure on the set of all partitions
which we shall call the {\bf Schur measure}. This measure depends
on countably many parameters.
In order to exhibit $\fM$ as a determinantal probability measure, as shown
by \ref b.Oko:wedge/, we first represent partitions as subsets of $\Z+1/2$
as follows.
To a partition $\la$ we associate a subset
$$
\fS(\la)=\{\la_i-i+1/2\}\subset \Z+{1 \over 2} \,.
$$
For example, $\fS(\emptyset)=\{-\frac 12,-\frac 32,-\frac 52\}$.
\msnote{???}
These are
known as the {\bf modified Frobenius coordinates of
$\la$}, see \ref b.KO/.
Then the law of $\fS(\la)$ when $\la$ has the law $\fM$ is $\P^K$, where
$K$ is defined by
the generating function
$$
\KK(z,w)=\sum_{i,j\in\Z+{1 \over 2}} z^i w^{-j}\, K(i,j)
\,,
$$
which is in turn defined by
$$
\KK(z,w) =  {\sqrt{zw} \over z-w} \,  {J(z) \over J(w)}\,.
$$
for $|w|<|z|$ and
$$
J(z):=
\prod_i {1-\overline{x_i}/z \over 1-x_i\,z } \,.
$$
}

Our main results are as follows (see \ref t.Q/).
We show that
for any positive contraction, $Q$, 
the measure $\P^Q$ has a very strong negative association property called
``conditional negative associations with external fields". Also, if
$Q_1$ and $Q_2$ are commuting positive contractions and $Q_1 \le Q_2$,
then $\P^{Q_1}$ is stochastically dominated by $\P^{Q_2}$.
These properties, especially the former, are powerful tools, comparable to the
well-known FKG inequalities representing positive association that hold
for many models of statistical mechanics.
For example, these results
are crucial to most of the results in \ref b.LS:dyn/.
%where
%we study the case that $E = \Z^d$ and $Q$ commutes with the action of
%$\Z^d$ on $\ell^2(\Z^d)$ (so that the measure $\P^Q$ is $\Z^d$-invariant).
We also show that $\P^Q$ has a trivial (full) tail $\sigma$-field for all $Q$;
this was shown independently by \ref b.ShiTak:II/ 
when the spectrum of $Q$ lies in $(0, 1)$. 
(Ironically, our proof is based on the case when the spectrum of $Q$ equals
$\{ 0, 1 \}$; see below.)
Another main result, \ref t.morris/, is
perhaps best appreciated at this point by giving two of its consequences.
The first consequence is a theorem of \ref b.Morris/, which says that on any
graph whatsoever, a.s.\
each tree in the wired uniform spanning forest is recurrent for simple random
walk (see \ref s.back/ for definitions).
A second corollary is the following.
Let $\T := \R/\Z$ be the unit circle equipped with unit Lebesgue measure. For
a measurable function $f : \T \to \C$ and an integer $n$, the {\bf Fourier
coefficient} of $f$ at $n$ is 
$$
\widehat f(n) := \int_\T f(t) e^{-2\pi i n t} \,dt
\,.
$$
For the meaning of ``Beurling-Malliavin density", see \ref d.BM/.

\procl c.not0
Let $A \subset \T$ be Lebesgue measurable with measure $|A|$.
Then there is a set $S$ of Beurling-Malliavin density $|A|$ in $\Z$ such that
if $f \in L^2(\T)$ vanishes a.e.\ on $A$ and $\widehat f$ vanishes off $S$,
then $f = 0$ a.e.
Indeed,
let $\P^A$ be the determinantal probability measure on $2^\Z$
corresponding to the Toeplitz matrix $(j, k) \mapsto \widehat{\I{A}}(k-j)$.
Then $\P^A$-a.e.\ $S \subset \Z$ has this property.
\endprocl

This can be compared to a theorem of \ref b.BTz/, according to which there is
a set $S \subset \Z$ of density at least $2^{-8} |A|$ 
such that if $f \in L^2(\T)$ and $\widehat f$ vanishes off $S$, then 
$$
|A|^{-1} \int_A |f(t)|^2 \,dt \ge 2^{-16} \|f\|_2^2
\,.
$$
Here, ``density'' is understood in the ordinary sense, but $S$ is found
so that the density of $S$ is also equal to its Schnirelman density.
It would be interesting to find a common strengthening of \ref c.not0/ and
the theorem of \refbauthor{BTz}.
Note that for \ref c.not0/, it is important to use a strong notion of density,
as otherwise a subset of $\N$ would be a trivial example for any $|A| \le 1/2$.
(\ref b.BTz/ find $S \subset \N$, but they could just as well have chosen $S$
in $\Z$ with the same density.)

As hinted at above, we shall base all our results on the case where $Q$ is an
orthogonal projection.
This is accomplished by a new and simple reduction method from positive
contractions to orthogonal projections.
A main tool in the case of orthogonal projections is a new representation of
the associated probability measure via exterior algebra.
We are also led naturally in this case to relations with matroids.
Finally, we detail a large number of new questions and conjectures
that we find quite
intriguing. These range over a number of areas of mathematics,
including 
Hilbert spaces, combinatorics, homology, and group representations,
{\it inter alia}. 

We now give an outline of the rest of the paper.
We begin in \ref s.matroid/ with the relationship of determinantal probability
measures to matroids. Matroids provide a useful organizing framework as well as
an inspirational viewpoint.
Prior research has contributed to a deep understanding of
the particular case of
the uniform spanning tree measure (\ref x.ust/)
and its extensions to infinite graphs.
Since generalizing this understanding is one goal of the present investigations
and of several open questions, we provide a quick summary of the relevant
facts for uniform spanning trees and forests in \ref s.back/.
%The easiest way to work with determinantal probabilities is to use the exterior
%algebra $\Ext\big(\ell^2(E)\big)$ of $\H$ with its usual inner product.
%In fact, \ref e.(3)/ \msnote{put this somewhere} drops out almost
%effortlessly, as does the Matrix-Tree Theorem and its generalization to
%regular matroids. In this way, we get the
%shortest route to many of the fundamental theorems connecting electrical
%networks and spanning trees, including Kirchhoff's theorem. 
For the sake of probabilists, 
we review exterior algebra in \ref s.ext/.
We use exterior algebra in \ref s.prmsr-d/ to give the basic
properties of determinantal probability measures and in \ref s.events/ to
prove our stochastic comparison inequalities.
%In an appendix, \ref s.ctg/, we use the present developments to give simple
%proofs of known results, such as the Matrix-Tree Theorem, on finding the
%cardinality of $\B$ in a regular matroid.
%
%We shall see that
%events for the probability space $(\B, \P^H)$ correspond to certain
%subspaces of $\Ext\big(\ell^2(E)\big)$. 
%The event is increasing iff the subspace is an
%ideal (with respect to exterior product). 
%Using this, we show in \ref s.events/ that
%if $H_1 \subset H_2$, then $\P^{H_1}$ is stochastically dominated by
%$\P^{H_2}$.
%This is a crucial tool in \ref b.LS:dyn/.
%It implies that for every probability measure $\P^H$, an increasing
%event $\ev A$ that ignores an element $e \in E$ is negatively correlated
%with the event that $\{e \in \ba\}$; this was proved by \ref b.FedMih/
%for the uniform measure on bases of ``balanced" matroids, a class of
%matroids containing the regular ones.
%We also use it to prove
%a more general negative association property.
%
The case of infinite dimensions is treated in \ref s.infinite/.
%, where we
%shall prove \ref t.basis/ stated above, as well as
%our extension of Morris's theorem alluded to above and
%the fact that every measure $\P^H$ has a trivial tail $\sigma$-field.
In \ref s.contract/, we explain how to associate probability measures to
positive contractions and show how results in this more
general setting follow easily from results for the special case of orthogonal
projections.
We outline many areas of open questions in Sections 
\briefref s.applic:thy/--\briefref s.applic:dyn/.

%We shall give a new simple proof of that fact in 
%Our new proof allows us easily to generalize results that we prove
%for the important special case where $Q$ is an orthogonal projection.

\bsection{Matroids}{s.matroid}

The set of spanning trees of a finite connected graph is not only the
best-known ground set for discrete determinantal probability measures,
but also the most
well-known example of the set of bases in a matroid.
Indeed, it was one of the two principal motivating examples for the
introduction of the theory of matroids by \ref b.Whitney:mat/.
Moreover, matroids are intrinsically linked to determinantal probability
measures.
To see this, recall the definition of a matroid.

A matroid is a simple combinatorial object satisfying just one axiom (for
more background, see \ref b.Welsh/ or \ref b.Oxley/).
Namely, let $E$ be a finite set %, called the {\bf ground set}, 
and let $\B$ be a nonempty collection of subsets of $E$.
We call the pair $\M :=
(E, \B)$ a {\bf matroid} with {\bf bases} $\B$ if the following exchange
property is satisfied:
$$
\forall B, B' \in \B \ \ \forall e \in B\setminus B' \ \ \exists e' \in
B'\setminus B \quad (B \setminus \{e\}) \cup \{e'\} \in \B \,.
$$
All bases have the same cardinality, called the {\bf rank} of the matroid.
Two fundamental examples are:
\beginbullets

$E$ is the set of edges of a finite connected graph $G$ and $\B$ is the set
of spanning trees of $G$.

$E$ is a set of vectors in a vector space and $\B$ is the collection of
maximal linearly independent subsets of $E$, where ``maximal'' means with
respect to inclusion.

\endbullets

The first of these two classes of examples is called {\bf graphical}, while
the second is called {\bf vectorial}.
The {\bf dual} of a matroid $\M = (E, \B)$ is the matroid $\M^\perp := (E,
\B^\perp)$, where $\B^\perp := \{E \setminu B \st B \in \B\}$. The dual matroid
is also called the {\bf orthogonal} matroid.

Any matroid that is isomorphic to a
vectorial matroid is called {\bf representable}. Graphical matroids are
{\bf regular}, i.e., representable over every field.

Each representation of a vectorial matroid over a subfield of the complex
numbers $\C$ gives rise to a determinantal probability measure in the
following fashion. 
Let $\M = (E, \B)$ be a matroid of rank $r$.
If $\M$ can be represented over a field $\Fd$,
then the usual way of specifying such a representation 
is by an $(s \times E)$-matrix $M$ whose columns are
the vectors in $\Fd^s$ representing $\M$ in the usual basis of $\Fd^s$.
For what follows below, in fact any basis of $\Fd^s$ may be used.
One calls $M$ a {\bf coordinatization matrix} of $\M$.
Now the rank of the matrix $M$ is also $r$.
%Let the columns of $M$ be $v_e$ ($e \in E$).
%The span of $\{v_e\}$ is $r$-dimensional.
If a basis for the column space is used instead of the usual basis of
$\Fd^s$, then we may take $s = r$.
In any case, the row space $H \subseteq \Fd^E$ of $M$ is $r$-dimensional and
we may assume that the first $r$ rows of $M$ span $H$.

When $\Fd \subseteq \C$, a determinantal probability measure $\P^H$
corresponding to the row space $H$ (or indeed to any subspace $H$ of
$\ell^2(E;\, \Fd)$) can be defined via any of several equivalent formulae.
Conceptually the simplest definition is to use \ref e.DPM/ with $Q$ being
the matrix of the orthogonal projection $P_H : \H \to H$, where the matrix of
$P_H$ is defined with respect to the usual basis of $\H$.
If, however, the coordinatization matrix $M$ is more available than $P_H$, then
one can proceed as follows.

For an $r$-element subset $B$ of $E$, let $M_B$ denote the $(r \times
r)$-matrix determined by the first $r$ rows of $M$ and those columns of $M$
indexed by $e \in B$.
Let $M_{(r)}$ denote the matrix formed by the {\it entire\/} first $r$ rows of
$M$.
In \ref s.prmsr-d/, we shall see that 
$$
\P^H[B] = |\det M_B|^2/\det(M_{(r)} M^*_{(r)})
\,, \label e.rowform
$$
where the superscript $*$ denotes adjoint.
(One way to see that this
defines a probability measure is to use the Cauchy-Binet formula.)
%Thus, $\P^H[B]$ is proportional to the square of the $r$-dimensional volume of
%the parallelepiped determined by $\Seq{v_e \st e \in B}$. 
%Use the form with $s = r$ to see this.
Identify each $e \in E$ with the corresponding unit basis vector of $\H$.
Equation \ref e.rowform/ shows that the matrix $M$ is a
coordinatization matrix of $\M$ iff the row space $H$ is a {\bf subspace
representation} of $\M$ in the sense that for $r$-element subsets $B
\subseteq E$, we have $B \in \B$ iff $P_H B$ is a basis for $H$. 
(We shall also say simply that $H$ {\bf represents} $\M$.)

Now use row operations to transform the first $r$ rows of $M$ to an
orthonormal basis for $H$.
Let $\widehat M$ be the $(r \times E)$-matrix
formed from these new first $r$ rows.
We shall see that 
$$
\widehat M^* \widehat M \hbox{ is the transpose of the matrix of } P_H
\,.
\label e.MM
$$

%We shall see that for any $A \subseteq E$, the analogue of \ref e.(3)/
%is
%$$
%\P^H[A \subseteq \ba] = \det [(\widehat M^* \widehat M)_{e, e'}]_{e, e' \in
%A}
%\,,
%\label e.rowTCT
%$$
%where the matrix on the right-hand side is the $(A \times A)$-minor of the
%matrix $\widehat M^* \widehat M$.

We now interpret all this for the case of a graphical matroid.
Given a finite connected
graph $G = (\vertex, \edge)$, choose an orientation
for each edge. Given $x \in \vertex$, define the {\bf star} at $x$ to be
the vector $\star_x := \star_x^G := \sum_{e \in \edge} a(x,e) e \in
\ell^2(\edge)$, where 
$$
a(x,e) := \cases{
0 &if $x$ is not incident to $e$,\cr
1 &if $x$ is the tail of $e$,\cr
-1 &if $x$ is the head of $e$\cr
}$$
are the entries of the $\vertex \times\edges$ {\bf vertex-edge incidence
matrix}.
Let $\STAR := \STAR(G)$ be the linear span of the stars, which is the same as
the row space of $[a(\cbuldot,\cbuldot)]$.
The standard coordinatization matrix of the matroid of spanning trees, the
graphic matroid of $G$, is $[a(\cbuldot,\cbuldot)]$. It is
easy to check that this does represent the graphic matroid.
It follows that $\STAR$ is a subspace representation of the graphic matroid,
which is also well known.
It is further well known that the matrix of $P_H$ is the transfer current
matrix $Y$; see, e.g., \BLPSusf.
It is well known that the numerator of the right-hand side of \ref
e.rowform/ is 1 when $B$ is a spanning tree of $G$ and 0 otherwise.
It follows that $\P^\STAR$ is the uniform measure on spanning trees of $G$.
It also follows that the denominator of the right-hand side of \ref e.rowform/
is the number of spanning trees.
To see what the matrix $M_{(r)} M^*_{(r)}$ is that appears
in the denominator of the right-hand side of \ref e.rowform/,
let the vertices of the graph be $v_0, \ldots, v_r$.
For $1 \le i, j \le r$ only, let $b_{i, j}$ be $-1$ if $v_i$ and $v_j$ are
adjacent, the degree of $v_i$ if $i = j$, and 0 otherwise.  This is the
same as the combinatorial Laplacian after the row and column corresponding
to $v_0$ are removed.
Then if the first $r$ rows of $M$ correspond to $v_i$ for $1 \le i \le r$, we
obtain that $M_{(r)} M^*_{(r)}$ is $[b_{i,j}]$.
The fact that the determinant of this matrix is the number of spanning trees
is known as the Matrix-Tree Theorem.
As we shall see in \ref c.uniform-reg/, those matroids $(E, \B)$ having a
subspace representation $H$ whose associated probability measure $\P^H$ is
the uniform measure on $\B$ are precisely the regular matroids in the case
of real scalars and complex unimodular matroids (also called
sixth-root-of-unity matroids) in the case of complex scalars.
%Several new avenues of investigation are suggested at the end of our paper
%that connect with various fields of mathematics. 
%Some of these are explored in \ref b.LS:dyn/.

\comment{
Every real matroid is orientable: Given a set of vectors $S$ in $\R^E$, let
the span of $S$ be $H$. Let $r := \dim H$. For $s_1, \ldots, s_r \in S$,
define $\chi(s_1, \ldots, s_r) := \sgn \ip{\bigwedge s_i, \xi_H}$.
}

There is a relationship to oriented matroids that may be worth noting.
In case the field of
a vectorial matroid is ordered, such as $\R$, there is an additional structure
making it an oriented matroid. Rather than define oriented matroids in
general, we give a definition in the representable case that is convenient
for our purposes: a {\bf real oriented matroid} $\chi$ of rank $r$ on $E =
\{1, 2, \ldots, n\}$ is a map $\chi : E^r \to \{-1, 0, +1\}$ such that for
some independent vectors $v_1, \ldots, v_r \in \R^E$,
$$
\chi(k_1, \ldots, k_r) = \sgn \det [\ip{v_{i}, e_{k_j}}]_{i, j \le r}\,, 
\label e.(2)
$$
where $e_1, \ldots, e_n$ is the standard basis of $\R^E$ and $\ip{\cbuldot,
\cbuldot}$ is the usual inner product in $\R^E$.  
The sets $\{k_1, \ldots, k_r\}$ with
$\chi(k_1, \ldots, k_r) \ne 0$ are the bases of a matroid on $E$.  
\par
Let $H$ be the linear span of the vectors $v_k$ appearing in the preceding
definition.  Choosing a different basis for $H$ and defining a new $\chi'$
in terms of this new basis by \ref e.(2)/ will lead only to $\chi' = \pm \chi$. 
Of course, the determinants
themselves in \ref e.(2)/ can change dramatically. However, if $\{v_1, \ldots,
v_r\}$ are orthonormal, then the determinants are fixed up to this sign
change and they give a determinantal probability measure.
Namely, 
if $\ba$ denotes a random base chosen with probability 
$$
\P^H[\ba = \{k_1, \ldots, k_r\}] 
  := (\det [\ip{v_{i}, e_{k_j}}]_{i, j \le r})^2\,,  \label e.detpr
$$
then this agrees with \ref e.rowform/, as shown in \ref s.prmsr-d/.

\bsection{Uniform Spanning Forests}{s.back}

The subject of random spanning trees of a graph goes back to \ref
b.Kirchhoff/,
who discovered several relations to electrical networks, one of
them mentioned already in \ref s.intro/.
His work also suggested the Matrix-Tree Theorem. % (\ref t.matrix-tree/).

In recent decades, computer scientists have developed various algorithms
for generating spanning trees of a finite graph at random according to the
uniform measure.  In particular, such algorithms are closely related to
generating states at random from any Markov chain. See \ref b.ProppWil/
for more on this issue. 

Early algorithms for generating a random spanning tree used the Matrix-Tree
Theorem. A much better algorithm than these early ones, especially for
probabilists, was introduced by \ref b.Aldous:ust/ and \ref b.Broder/. It says
that if you start a simple random walk at {\it any\/} vertex of a finite
graph $G$ and draw every edge it crosses except when it would complete a
cycle (i.e., except when it arrives at a previously-visited vertex), then
when no more edges can be added without creating a cycle, what will be
drawn is a uniformly chosen spanning tree of $G$. This beautiful algorithm
is quite efficient and useful for theoretical analysis, yet \ref
b.Wilson:gen/
found an even better one. It is lengthier to describe, so we forego
its description, but it is faster and more flexible.

The study of the analogue on an infinite graph of a uniform spanning tree
was begun by \ref b.Pemantle:ust/ at the suggestion of the present author.
Pemantle showed
that if an infinite connected graph $G$ is exhausted by a sequence of finite
connected subgraphs $G_n$, then the weak limit of the uniform spanning tree
measures on $G_n$ exists. 
However, it may happen that the limit measure is not supported on trees,
but on forests.
This limit measure is now called the {\bf free uniform spanning forest} on
$G$, denoted $\fusf$. Considerations of electrical networks play the
dominant role in the proof of existence of the limit.
If $G$ is itself a tree, then this measure is
trivial, namely, it is concentrated on $\{G\}$. Therefore, \ref b.Hag:rcust/
introduced another limit that had been considered more implicitly by \ref
b.Pemantle:ust/ on $\Z^d$, namely, the weak limit of the uniform
spanning tree measures on $G_n^*$, where $G_n^*$ is the graph $G_n$ with
its boundary identified to a single vertex. As \ref b.Pemantle:ust/ showed,
this limit also always exists
on any graph and is now called the {\bf wired uniform spanning forest},
denoted $\wusf$.  

In many cases, the free and the wired limits are the same. In
particular, this is the case on all euclidean lattices such as $\Z^d$.
The general question of when the free and wired uniform spanning forest
measures are the same turns out to be quite interesting: The measures are the
same iff there are no nonconstant harmonic Dirichlet functions on $G$ (see
\BLPSusf).

Among the notable results that \ref b.Pemantle:ust/ proved was that on $\Z^d$,
the uniform spanning forest is a single tree with a single end a.s.\ if $2
\le d \le 4$; when $d \ge 5$, there are infinitely many trees a.s., each
with at most two ends. \BLPSusf\ shows that, in fact, not only in
$\Z^d$, but on any Cayley graph of a finitely generated group other than a
finite extension of $\Z$, each tree in the $\wusf$ has but a single end a.s.
One of Pemantle's main tools was the Aldous/Broder algorithm, while \BLPSusf\
relied on a modification of Wilson's algorithm.

Another result of \ref b.Pemantle:ust/ was that on $\Z^d$, the tail of the
uniform spanning forest measure is trivial. This is extended to both the
free and the wired spanning forests on all graphs in \BLPSusf.

In the paper \BLPSusf, it was noted that the Transfer Current Theorem
extends to the free and wired spanning forests if one uses the free and
wired currents, respectively. 
To explain this, note that the orthocomplement of the row space $\STAR(G)$ of
the vertex-edge incidence matrix of a finite graph $G$ is the kernel, denoted
$\CYCLE(G)$, of the matrix.
We call $\STAR(G)$ the {\bf star space} of $G$ and $\CYCLE(G)$ the {\bf cycle
space} of $G$.
For an infinite graph $G = (\vertex, \edges)$ exhausted by finite subgraphs
$G_n$, we let $\STAR(G)$ be the closure of $\bigcup \STAR(G_n^*)$ and
$\CYCLE(G)$ be the closure of $\bigcup \CYCLE(G_n)$, where we take the closure
in $\ell^2(\edges)$.
Then $\wsf$ is the determinantal probability measure corresponding to the
projection $P_{\STAR(G)}$, while $\fsf$ is the determinantal probability measure
corresponding to $P_{\CYCLE(G)}^\perp := P_{\CYCLE^\perp}$.
In particular, $\wsf = \fsf$ iff $\STAR(G) = \CYCLE(G)^\perp$.

While the wired spanning forest is quite well understood,
the free spanning forest measure is in general poorly understood. 
A more detailed summary of uniform spanning forest measures can be found in
\ref b.Lyons:bird/.

\bsection{Review of Exterior Algebra}{s.ext}

Let $E$ be a finite or countable set.  Identify $E$ with the standard
orthonormal basis of the real or complex Hilbert space $\ell^2(E)$.
For $k \ge 1$, let $E_k$ denote a collection
of ordered $k$-element subsets of $E$ such that
each $k$-element subset of $E$ appears exactly once in $E_k$ in some
ordering. Define 
$$
\Lambda^k E := \bigwedge\nolimits^{\!k} \H
:= \ell^2\Bigl( \big\{e_1 \wedge \cdots
\wedge e_k \st \Seq{e_1, \ldots, e_k} \in E_k \big\} \Bigr)
\,.
$$
If $k > |E|$, then $E_k = \emptyset$ and $\Lambda^k
E = \{0\}$.  We also define $\Lambda^0 E$ to be 
the scalar field, $\R$ or $\C$.
The elements of $\Lambda^k E$ are called {\bf
multivectors} of {\bf rank} $k$, or {\bf $k$-vectors} for short. We then
define the {\bf exterior \rm (or \bf wedge\rm) \bf product}
of multivectors in the usual alternating
multilinear way: $\bigwedge_{i=1}^k e_{\sigma(i)} = \sgn(\sigma)
\bigwedge_{i=1}^k e_i$ for any permutation $\sigma$ of $\{1, 2, \ldots,
k\}$ and $\bigwedge_{i=1}^k
\sum_{e\in E'} a_i(e) e = \sum_{e_1, \ldots, e_k \in E'} \prod_{j=1}^k
a_j(e_j) \bigwedge_{i=1}^k e_i$ for any scalars $a_i(e)$ ($i \in [1, k],\;
e \in E'$), and any finite $E' \subseteq E$.
(Thus, $\bigwedge_{i=1}^k e_i = 0$ unless
all $e_i$ are distinct.) 
The inner product on $\Lambda^k E$ satisfies
$$
\ip{u_1 \wedge \cdots \wedge u_k, v_1 \wedge \cdots \wedge v_k}
=
\det\big[\ip{u_i, v_j}\big]_{i, j \in [1, k]}
\label e.ipdet
$$
when $u_i$ and $v_j$ are 1-vectors.
(This also shows that the inner product on $\Lambda^k E$ does not
depend on the choice of orthonormal basis of $\H$.)
We then define the {\bf exterior} (or {\bf Grassmann}) {\bf algebra}
$\Ext\big(\H\big) := \Ext(E) := \bigoplus_{k \ge 0} \Lambda^k E$, where the
summands are declared orthogonal, making it into a Hilbert space.
(Throughout the paper, $\oplus$ is used to indicate the sum
of orthogonal summands, or, if there are an infinite number of orthogonal
summands, the closure of their sum.)
Vectors $u_1, \ldots, u_k \in \H$ are linearly independent iff $u_1
\wedge \cdots \wedge u_k \ne 0$.
For a $k$-element subset $A \subseteq E$ with ordering $\Seq{e_i}$ in $E_k$,
write
$$
\theta_A := \bigwedge_{i=1}^k e_i
\,.
$$
We also write 
$$
\bigwedge_{e \in A} f(e) := \bigwedge_{i=1}^k f(e_i)
$$
for any function $f : E \to \H$.

If $H$ is a (closed) linear subspace of $\H$, then we identify $\Ext(H)$
with its inclusion in $\Ext(E)$. That is, $\bigwedge^k H$ is the closure of
the linear span of $\{v_1 \wedge \cdots \wedge v_k \st v_1, \ldots, v_k \in
H\}$.
In particular, if $\dim H = r < \infty$, then
$\bigwedge^r H$ is a 1-dimensional subspace of $\Ext(E)$; denote by
$\xi_H$ a unit multivector in this subspace, which is unique up to sign in
the real case and unique up to a scalar unit-modulus factor in the complex
case, i.e., up to {\bf signum}.
We denote by $P_H$ the orthogonal projection onto $H$ for any subspace $H
\subseteq \H$ or, more generally, $H \subseteq \Ext(E)$. 

\procl l.projection
For any subspace $H \subseteq \H$, any $k \ge 1$, and any $u_1, \ldots, u_k
\in \H$,
$$
P_{\Ext(H)} (u_1 \wedge \cdots \wedge u_k) = (P_H u_1) \wedge \cdots
\wedge (P_H u_k)
\,.
$$
\endprocl

\proof 
Write 
$$
u_1 \wedge \cdots \wedge u_k
=
(P_H u_1 + P^\perp_H u_1) \wedge \cdots \wedge (P_H u_k + P^\perp_H u_k)
$$
and expand the product. All terms but $P_H u_1 \wedge \cdots \wedge P_H
u_k$ have a factor of $P^\perp_H u$ in them, making them orthogonal to
$\Ext(H)$ by \ref e.ipdet/. 
\Qed

A multivector is called {\bf simple} or {\bf decomposable} if it is the
wedge product of 1-vectors.  \ref b.Whitney:book/, p.~49, shows that 
$$
\|\u \wedge \v \| \le \|\u\| \|\v\| \quad\hbox{ if either $\u$ or $\v$ is
simple}. \label e.whitney
$$
%This implies Hadamard's inequality for determinants:
%for $n$ vectors $u_1, u_2, \ldots, u_n$ belonging to $n$-dimensional space
%$\C^n$, the matrix $A$ whose $i$th row is $u_i$ has (by \ref e.ipdet/ and \ref
%e.whitney/)
%$$
%|\det A|
%=
%\|u_1 \wedge u_2 \wedge \cdots \wedge u_n\| 
%\le
%\prod_{i=1}^n \|u_i\|
%\,.
%\label e.hadamard
%$$

We shall use the {\bf interior product} defined by duality:
$$
\ip{\u \vee \v, \w} = \ip{\u, \w \wedge \v}
\qquad (\u \in \Lambda^{k+l} E ,\ \v \in \Lambda^l E, \ \w \in
\Lambda^k E)\,.
$$
In particular, if $e\in E$ and $\u$ is a multivector that does not contain
any term with $e$ in it (that is, $\u\in\Ext(e^\perp)$), then $(\u \wedge
e) \vee e = \u$ and $\u \vee e = 0$.
More generally, if $v \in \H$ with $\|v\|=1$ and $\u \in \Ext(v^\perp)$,
then $(\u \wedge v) \vee v = \u$ and $\u \vee v = 0$.
Note that the interior product is sesquilinear, not bilinear, over $\C$.

For $e \in E$, write $[e]$ for the subspace of scalar multiples of $e$ in $\H$.
If $H$ is a finite-dimensional subspace of $\H$ and $e \notin H$, then 
$$
\xi_H \wedge e = \|P_H^\perp e\|\;\xi_{H + [e]}
\label e.Hwedge
$$
(up to signum).
To see this, let $u_1, u_2, \ldots, u_r$ be an orthonormal basis of $H$,
where $r = \dim H$.
Put $v := P_H^\perp e /\|P_H^\perp e\|$.
Then $u_1, \ldots, u_r, v$ is an orthonormal basis of $H + [e]$, whence
$$
\xi_{H + [e]} 
=
u_1 \wedge u_2 \wedge \cdots \wedge u_r \wedge v
=
\xi_H \wedge v
=
\xi_H \wedge e/\|P_H^\perp e\|
$$
since $\xi_H \wedge P_H e = 0$.
This shows \ref e.Hwedge/.
Similarly, if $e \notin H^\perp$, then 
$$
\xi_H \vee e = \|P_H e\|\;\xi_{H \cap e^\perp}
\label e.Hvee
$$
(up to signum).
Indeed, put $w_1 := P_H e /\|P_H e\|$. Let $w_1, w_2, \ldots, w_r$ be an
orthonormal basis of $H$ with $\xi_H = w_1 \wedge w_2 \wedge \cdots \wedge
w_r$.
Then 
$$
\xi_H \vee e
=
\xi_H \vee P_H e
=
(-1)^{r-1} \|P_H e\|\; w_2 \wedge w_3 \cdots \wedge w_r
$$
(up to signum), as desired.

Finally, we claim that 
$$
\all {u, v \in \H}\quad \bigip{\xi_H \vee u, \xi_H \vee v} 
=
\ip{P_H v, u}
\,.
\label e.reverse
$$
Indeed, $\xi_H \vee u = \xi_H \vee P_H u$, so this is equivalent to
$$
\bigip{\xi_H \vee P_H u, \xi_H \vee P_H v} 
=
\ip{P_H v, P_H u}
\,.
$$
Thus, it suffices to show that
$$
\all {u, v \in H}\quad \bigip{\xi_H \vee u, \xi_H \vee v} 
=
\ip{v, u}
\,.
$$
By sesquilinearity, it suffices to show this for $u, v$ members of an
orthonormal basis of $H$.
But then it is obvious.

\comment{The following is not used.
If $E$ is finite, then the {\bf dual} $\D \u$ of a multivector $\u \in
\Ext(E)$ is
$\D \u = \xi_{\H} \vee \u$.
For example, if $H_2 = H_1^\perp$, then
$\D \xi_{H_1} = \xi_{\H} \vee \xi_{H_1} = (\pm \xi_{H_1} \wedge
\xi_{H_2}) \vee \xi_{H_1} = \pm \xi_{H_2}$.
}

\comment{Here is how to get $\det(AB) = \det A \cdot \det B$: First, that
$\det A = \det A^*$ follows from symmetry of inner product. Now let the
rows of $A$ be $v_i$ and the columns of $B$ be $w_i$. Let $\theta :=
\theta_E$.  Then $\det (AB) = \det \ip{v_i, w_j} = \ip{\bigwedge v_i, \bigwedge
w_i}$. Substitute $\bigwedge v_i = \ip{\bigwedge v_i, \theta} \theta$  and the
same for $\bigwedge w_i$. We obtain $\det (AB) = \det A \det B^*$.}

For a more detailed presentation of exterior algebra, see \ref
b.Whitney:book/.

\bsection{Probability Measures on Bases}{s.prmsr-d}

Any unit vector $v$ in a Hilbert space with orthonormal basis $E$ gives a
probability measure $\P^v$ on $E$, namely, $\P^v[\{e\}] := |\ip{v, e}|^2$
for $e \in E$. Applying this simple idea to multivectors instead, we shall
obtain the probability measures $\P^H$ of \ref s.intro/.

Let $H$ be a subspace representation of
the matroid $\M = (E, \B)$ of rank $r$.
Then the non-0 coefficients in $\xi_H$
with respect to the standard basis of $\Lambda^r E$ are exactly those of
the multivectors $\theta_B$ ($B \in \B$).
Indeed, by \ref l.projection/,
the coefficient in $\xi_H$ of $\theta_B = \bigwedge_{i=1}^r e_i$ satisfies
$$
\left|\leftip{\xi_H, \theta_B}\right|
=
\|P_{\Ext(H)} (\bigwedge_{i=1}^r e_i) \| 
= 
\| \bigwedge_i P_H e_i \|\,,
$$
which is non-0 iff $\Seq{P_H e_i}$ are linearly independent.
Furthermore,
we may define the probability measure $\P^H$ on $\B$ by
$$
\P^H[B] 
:=
|\leftip{\xi_H, \theta_B}|^2
\,.  \label e.xiHpr
$$
%$$
%\P^H[\{e_1, \ldots, e_r\}] 
%:=
%\left|\leftip{\xi_H, \bigwedge_{i=1}^r e_i}\right|^2
%\,.  \label e.xiHpr
%$$
We may also write \ref e.xiHpr/ as
$$
\xi_H = \sum_{B \in \B} \epsilon_B \sqrt{\P^H[B]} \theta_B
$$
for some $\epsilon_B$ of absolute value 1, or alternatively as
$$
\P^H[\ba= B]
=
\|P_{\Ext(H)} \theta_{B}\|^2
=
\ip{P_{\Ext(H)} \theta_{B}, \theta_{B}}
\,.
$$

To see that \ref e.xiHpr/ agrees with \ref e.detpr/, we just use \ref
e.ipdet/ and the fact that $\xi_H = c \bigwedge_i v_i$ for some scalar $c$
with $|c| = 1$.

To show that the definition \ref e.xiHpr/ satisfies \ref e.DPM/ for the
matrix of $P_H$, observe that
$$
\P^H[\ba = B]
=
\bigip{P_{\Ext(H)} \theta_B, \theta_B}
=
\leftip{\bigwedge_{e \in B} P_H e, \bigwedge_{e \in B} e}
=
\det [\ip{P_H e, f}]_{e, f \in B}
$$
by \ref e.ipdet/.
This shows that \ref e.DPM/ holds for $A \in \B$ since $|\ba| = r$ $\P^H$-a.s.
We now prove the general case by proving an extension of it.
We shall use the convention that $\theta_ \emptyset := 1$ and $\u \wedge 1 :=
\u$ for any multivector $\u$.

\procl t.genprs
Suppose that $A_1$ and  $A_2$ are (possibly empty) subsets of a finite set $E$. 
We have
$$
\P^H[A_1 \subseteq \ba, A_2 \cap \ba = \emptyset]
=
\bigip{P_{\Ext(H)} \theta_{A_1} \wedge P_{\Ext(H^\perp)} \theta_{A_2},
\theta_{A_1} \wedge \theta_{A_2}}
\,.
\label e.genprs
$$
In particular, for any $A \subseteq E$, we have
$$
\P^H[A \subseteq \ba]
=
\|P_{\Ext(H)} \theta_A\|^2
\,.
\label e.included
$$
\endprocl

\procl r.gentct
%The Transfer Current Theorem is the special case
%where $H = \STAR$ for a graphic matroid.
The property \ref e.genprs/ (in an equivalent form) appears in several
places in the literature, including Theorem 3.1 of \ref b.ShiTak:ann/ and
Proposition A.8 of \ref b.BOO:plancherel/.
Usually, it is derived from a different definition of the measure $\P^H$.
For the uniform measure on bases of regular matroids, the case of \ref
e.included/ where $|A| = 2$ is Theorem 2.1 of \ref b.FedMih/.
The general case is related to Theorem 5.2.1 of \ref b.Mehta:book/.
\endprocl

\proof
Both sides of \ref e.genprs/ are clearly 0 unless $A_1 \cap A_2 =
\emptyset$, so we assume this disjointness from now on.
Let $r := \dim H$.

Consider next the case where $A_1 \cup A_2 = E$. The left-hand side
is 0 unless $|A_1| = r$ since every base has $r$
elements. Also the right-hand side is 0 except in such a case since
$\Ext(H)$ has multivectors only of rank at most $r$ and
$\Ext(H^\perp)$ has multivectors only of rank at most $|E|-r$.
It remains to establish equality when $|A_1| = r$.
In that case, we have
\begineqalno
P_{\Ext(H)} \theta_{A_1} \wedge \theta_{A_2}
&=
P_{\Ext(H)} \theta_{A_1} \wedge 
\bigwedge_{e \in A_2} (P_H e + P^\perp_H e)
\cr&=
P_{\Ext(H)} \theta_{A_1} \wedge P_{\Ext(H^\perp)} \theta_{A_2}
\,,\cr
\endeqalno
as we see by expanding the product, using the fact that $\Ext(H)$ has
multivectors only of rank at most $r$, and using \ref l.projection/.
Therefore, 
\begineqalno
\ip{P_{\Ext(H)} \theta_{A_1} \wedge P_{\Ext(H^\perp)} \theta_{A_2},
\theta_{A_1} \wedge \theta_{A_2}}
&=
\ip{P_{\Ext(H)} \theta_{A_1} \wedge \theta_{A_2},
\theta_{A_1} \wedge \theta_{A_2}}
\cr&=
\ip{P_{\Ext(H)} \theta_{A_1}, \theta_{A_1}}
\cr&=
\P^H[\ba= A_1]
\cr&=
\P^H[A_1 \subseteq \ba, A_2 \cap \ba = \emptyset]
\,.\cr
\endeqalno

This establishes the result when $E$ is partitioned as $A_1 \cup A_2$.
In the general case, we can decompose over partitions to write 
$$
\P^H[A_1 \subseteq \ba, A_2 \cap \ba = \emptyset]
=
\sum_{C \subseteq E \setminus (A_1 \cup A_2)}
\P^H[A_1 \cup C \subseteq \ba, 
\big(A_2 \cup (E \setminus C)\big) \cap \ba = \emptyset]
\qquad
\label e.1
$$
and we can also write 
\begineqalno
\biggip{P_{\Ext(H)} \theta_{A_1} \wedge P_{\Ext(H^\perp)} \theta_{A_2},
\theta_{A_1} \wedge \theta_{A_2}}
\hskip-2in&\hskip2in
=
\biggip{\bigwedge_{e \in A_1} P_{H} e \wedge 
\bigwedge_{e \in A_2} P_{H}^\perp e,
\theta_{A_1} \wedge \theta_{A_2}}
\cr&=
\biggip{\bigwedge_{e \in A_1} P_{H} e \wedge 
\bigwedge_{e \in A_2} P_{H}^\perp e \wedge
\bigwedge_{e \notin A_1 \cup A_2} e,
\theta_{A_1} \wedge \theta_{A_2} \wedge
\bigwedge_{e \notin A_1 \cup A_2} e}
\cr&=
\biggip{\bigwedge_{e \in A_1} P_{H} e \wedge 
\bigwedge_{e \in A_2} P_{H}^\perp e \wedge
\bigwedge_{e \notin A_1 \cup A_2} (P_H e + P_H^\perp e),
\theta_{A_1} \wedge \theta_{A_2} \wedge
\bigwedge_{e \notin A_1 \cup A_2} e}
\cr&=
\sum_{C \subseteq E \setminus (A_1 \cup A_2)}
\biggip{\bigwedge_{e \in A_1 \cup C} P_{H} e \wedge 
\bigwedge_{e \in A_2 \cup (E \setminus C)} P_{H}^\perp e,
\theta_{A_1 \cup C} \wedge \theta_{A_2 \cup (E \setminus C)}}
\,.
\label e.2
\cr
\endeqalno
For each $C$, the summand on the right of \ref e.1/ is equal to that on the
right of \ref e.2/ by what we have shown, whence the general case follows.
\Qed

We see immediately the relationship of orthogonality to duality: 

\procl c.dualrep
If a subspace $H \subseteq \H$ represents a matroid $\M$, then its orthogonal
complement $H^\perp$ represents the dual matroid $\M^\perp$: 
$$
\all {B \in 2^E}\quad \P^{H^\perp}[E \setminus B] = \P^H[B]
\,.
\label e.dualrep
$$
\endprocl

We now verify \ref e.rowform/ and \ref e.MM/.
Resuming the notation used there, we let the $i$th row of $M$ be $m_i$.
For some constant $c$, we thus have
$$
\xi_H = c \bigwedge_{i=1}^r m_i
\,,
\label e.xi-M
$$
whence by \ref e.ipdet/,
$$
\P^H[B]
=
|\ip{\xi_H, \theta_B}|^2
=
|c|^2 
\left|\det\left[\ip{m_i, e}\right]_{i \le r,\ e \in B}\right|^2
=
|c|^2 |\det M_B|^2
\,.
$$
To complete the derivation of \ref e.rowform/, we must calculate $|c|^2$. 
For this, note that
$$
1 
=
\|\xi_H\|^2
=
|c|^2 \big\| \bigwedge_{i=1}^r m_i \big\|^2
=
|c|^2 \det\left[\ip{m_i, m_j} \right]_{i, j \le r}
=
|c|^2 \det(M_{(r)} M_{(r)}^*)
\,.
$$

We record explicitly the following consequence noted in \ref s.intro/:

\procl c.row-representation
A matrix $M$ is a coordinatization matrix of a matroid $\M = (E, \B)$ iff
the row space $H$ of $M$ is a subspace representation of $\M$.
\endprocl

This is also obvious from the following equivalences that hold 
for any scalars $a_e$:
\begineqalno
\all i \sum_e a_e \ip{m_i, e} = 0 
&\iff
\sum_e a_e e \in \ker M = H^\perp
\iff
P_H \sum_e a_e e = 0
\cr&\iff \sum_e a_e P_H e = 0
\,.
\endeqalno

By \ref c.row-representation/, it follows that
$\STAR$ is a subspace representation of the graphic matroid, as mentioned in
\ref s.intro/.

In order to demonstrate \ref e.MM/, first note that row operations on
$M$ do not change the row space, $H$. Now let $\Seq{w_i \st 1 \le i \le r}
\in \H$ be the rows of $\widehat M$ and let the entries be $\widehat M_{i,
e} = \ip{w_i, e}$. 
Then
$$
P_H e = 
\sum_{i=1}^r \ip{e, w_i} w_i
=
\sum_{i=1}^r \overline{\widehat M_{i, e}} w_i
\,,
$$
whence $\ip{P_H e, e'} = \sum_{i=1}^r \overline{\widehat M_{i, e}} \ip{w_i,
e'} = \sum_{i=1}^r \overline{\widehat M_{i, e}} \widehat M_{i, e'} =
(\widehat M^* \widehat M)_{e, e'}$. 
Since $\ip{P_H e, e'}$ is the $(e', e)$-entry of the matrix of $P_H$, we
obtain \ref e.MM/.

We call a matroid {\bf complex unimodular} if it has a coordinatization
matrix $M$ over $\C$ with $|\det M_B| = 1$ for all bases $B$. These are
known to be the same as the {\bf sixth-root-of-unity}
matroids, where the condition now is that it has a coordinatization matrix
$M$ over $\C$ with $(\det M_B)^6 = 1$ for all bases $B$ (see \ref b.COSW/
for this and other characterizations of such matroids).

\procl c.uniform-reg
A matroid $\M = (E, \B)$ has a representation via a subspace $H \subseteq
\ell^2(E;\,\R)$ with $\P^H$ equal to uniform measure on $\B$ iff $\M$ is a
regular matroid.
A matroid $\M = (E, \B)$ has a representation via a subspace $H \subseteq
\ell^2(E;\,\C)$ with $\P^H$ equal to uniform measure on $\B$ iff $\M$ is a
complex unimodular matroid.
\endprocl

\proof
Having a subspace representation that gives uniform measure is equivalent
to having a coordinatization matrix $M$ with the above matrices $M_B$
having determinants equal in absolute value for $B \in \B$. In case the
field is
$\R$, this is equivalent to requiring all $\det M_B = \pm1$.
From the above, we see that this is the same as having some $H$
representing $\B$ and some constant $c$ with $c \xi_H$ having all its
coefficients in $\{ 0, \pm1 \}$.
This condition is equivalent to regularity: see (3) and (4) of
Theorem 3.1.1 of \ref b.White:comb/.

The statement for complex representations is clear from the definition.
\Qed

\procl r.weight
In dealing with electrical networks that have general conductances (not
necessarily all 1), one uses not the uniform spanning tree measure, but the
weighted spanning tree measure, where the probability of a tree is
proportional to the product of the conductances of its edges. 
If the conductance, or weight, of the edge $e$ is $w(e)$, then one defines the
star at a vertex $x$ by $\star_x := \sum_{e \in \edge} \sqrt{w(e)} a(x, e) e$
and one defines $\STAR$ and $\CYCLE$ accordingly. The weighted spanning tree
measure is then $\P^\STAR$.
In fact,
given any positive weight function $w$ on a ground set $E$ and given any
measure $\P^H$ on a matroid, there is another subspace $H_w$ such that for
all bases $B$, we have $\P^{H_w}[B] = Z^{-1} \P^H[B] \prod_{e \in B} w(e)$,
where $Z$ is a normalizing constant. To see this, define the linear
transformation $D_w$ on $\H$ by $D_w(e) = \sqrt{w(e)} e$. Then it is not
hard to verify that the image $H_w$ of $H$ under $D_w$ has this property.
Similarly, if the matroid is represented by the matrix $M$ with row space
$H$, then the matrix $M D_w$ represents the same matroid but with the
transformed measure $\P^{H_w}$.
\endprocl

\bsection{Stochastic Domination and Conditioning}{s.events}

Let $2^E$ denote the collection of all subsets of $E$.
Any probability measure $\P$ on $\B$ extends to a probability measure on
the Borel $\sigma$-field of $2^E$ by setting $\P(\ev A):=\P(\ev A \cap
\B)$.
In this section, we compare different probability measures on $2^E$ with
$E$ finite. First, for $E$ finite,
to each event $\ev A \subseteq 2^E$, we associate the subspace
$$
S_{\ev A} := \hbox{the linear span of } \big\{\theta_C \st C \in \ev A\big\}
\subseteq \Ext(E)\,.
$$
Note that $C$ need not be in $\B$ here.
An event $\ev A$ is called {\bf increasing} if whenever $A \in \ev A$ and
$e \in E$, we have also $A \cup \{e\} \in \ev A$. Thus, $\ev A$ is
increasing iff the subspace $S_{\ev A}$ is an ideal (with respect to
exterior product). 

The probability measure $\P^H$ clearly satisfies
$$
\P^H(\ev A) = \|P_{S_{\ev A}} \xi_H\|^2
$$
for any event $\ev A$.

\procl l.ideal
Let $S \subseteq \Ext(E)$ be an ideal, $\u = \bigwedge_{i=1}^k u_i$,
and $\v = \bigwedge_{j=1}^l v_j$. If $\ip{u_i, v_j} = 0$ for all $i$ and $j$,
then
$$
\bigip{ (P_S \u) \wedge \v, P_S(\u \wedge \v) }
=
\|P_S \u\|^2 \|\v\|^2
\,.
$$
\endprocl

\proof
Since $(P_S \u) \wedge \v \in S$, we have
\begineqalno
\bigip{ (P_S \u) \wedge \v, P_S(\u \wedge \v) }
&=
\leftip{ P_S\big((P_S \u) \wedge \v\big), \u \wedge \v }
\cr&=
\bigip{ (P_S \u) \wedge \v, \u \wedge \v }
\cr&=
\ip{P_S \u, \u} \ip{\v, \v}
\cr&=
\|P_S \u\|^2 \|\v\|^2
\,. \Qed \cr
\endeqalno

Given two probability measures $\P^1$, $\P^2$ on $2^E$, we say that {\bf
$\P^2$ stochastically dominates $\P^1$} and write $\P^1 \preccurlyeq \P^2$ if
for all increasing events $\ev A$, we have $\P^1(\ev A) \le \P^2(\ev A)$.

\procl t.dominate
Let $E$ be finite and let
$H_1 \subset H_2$ be subspaces of $\H$. Then $\P^{H_1} \preccurlyeq
\P^{H_2}$.
\endprocl

\proof
Let $\ev A$ be an increasing event. Take an orthonormal basis $\Seq{u_i
\st i \le r_2}$ of $H_2$ such that $\Seq{u_i \st i \le r_1}$ is an orthonormal
basis of $H_1$ and such that $\xi_{H_j} = \bigwedge_{i=1}^{r_j} u_i$ for
$j=1, 2$. Apply \ref l.ideal/ to $\u := \bigwedge_{i=1}^{r_1} u_i$ and $\v
:= \bigwedge_{i=r_1+1}^{r_2} u_i$ to conclude that 
\begineqalno
\|P_{S_\A} \xi_{H_1}\|^2
&=
\|P_{S_\A} \xi_{H_1}\|^2 \|\v\|^2
=
\bigip{(P_{S_\A} \xi_{H_1}) \wedge \v, P_{S_\A} \xi_{H_2}}
\cr&\le
\|P_{S_\A} \xi_{H_1} \wedge \v\| \|P_{S_\A} \xi_{H_2}\|
\le
\|P_{S_\A} \xi_{H_1}\| \|\v\| \|P_{S_\A} \xi_{H_2}\|
\cr&=
\|P_{S_\A} \xi_{H_1}\| \|P_{S_\A} \xi_{H_2}\|
\cr
\endeqalno
by \ref e.whitney/, whence $\|P_{S_\A} \xi_{H_1}\| \le \|P_{S_\A}
\xi_{H_2}\|$, i.e., $\P^{H_1}(\A) \le \P^{H_2}(\A)$.  \Qed

A {\bf minor} of a matroid $\M = (E, \B)$ is one obtained from $\M$ by
repeated contraction and deletion: Given $e \in E$, the {\bf contraction}
of $\M$ along $e$ is $\M/e := (E, \B/e)$, where $\B/e := \{ B \in \B \st e
\in B\}$, while the {\bf deletion} of $\M$ along $e$ is $\M\setminu e := (E,
\B\setminu e)$, where $\B\setminu e := \{ B \in \B \st e \notin B\}$. Note
that, contrary to usual convention, we are keeping the same ground set,
$E$.

For any $F \subseteq E$, let $[F]$ be the closure of the linear span of the
unit vectors $\{e \st e \in F\}$.
%In particular, we shall write $[e]$ when $F = \{ e \}$.
We shall also write $[u]$ for the subspace spanned by any vector $u$.

\procl p.minors
Let $E$ be finite and $H$ be a subspace of $\H$.
For any $e \in E$ with $0 < \P^H[e \in \ba] < 1$, we have 
$$
\P^H(\; \cbuldot \mid e \in \ba) = \P^{(H \cap e^\perp)+[e]}(\; \cbuldot\; )
$$
and 
$$
\P^H(\;
\cbuldot \mid e \notin \ba) = \P^{(H + [e]) \cap e^\perp}(\; \cbuldot\; )
\,.
$$
In particular, $(H \cap e^\perp)+[e]$ represents $\M/e$ and $(H + [e])
\cap e^\perp$ represents $\M\setminu e$.
Moreover, signs (or, in the complex case, unit scalar factors)
can be chosen for $\xi_{(H \cap e^\perp)+[e]}$ and
$\xi_{(H + [e]) \cap e^\perp}$ so that for all $B \in \B$,
$$
\bigip{\theta_B, \xi_{(H \cap e^\perp)+[e]}} = \cases{
   \ip{\theta_B, \xi_H}/\| P_H e \|  &if $e \in B$,\cr
   0                                 &if $e \notin B$\cr}
\label e.mincoeff1
$$
and
$$
\bigip{\theta_B, \xi_{(H + [e]) \cap e^\perp}} = \cases{
   \ip{\theta_B, \xi_H}/\| P_H^\perp e \|  &if $e \notin B$,\cr
   0                                 &if $e \in B$.\cr}
\label e.mincoeff2
$$
\endprocl

\proof
Clearly $\ip{\theta_B, \xi_H} = \bigip{\theta_B, (\xi_H \vee e) \wedge e}$ when
$e \in B$, while the right-hand side of this equation is 0 for $e \notin B$. 
Applying \ref e.Hvee/ followed by \ref e.Hwedge/, we obtain \ref e.mincoeff1/.
The proof of \ref e.mincoeff2/ is similar.
These imply the formulas for the conditional probabilities.
\Qed

Given a disjoint pair of possibly infinite sets $A, B \subseteq E$, define 
$$
H_{A, B} := \big( ( H \cap A^\perp ) + [A \cup B] \big) \cap B^\perp
\,.
\label e.HAB
$$
It is straightforward to verify that 
$$
H_{A, B} = \big( ( H + [B] ) \cap (A \cup B)^\perp \big) + [A]
\,.
\label e.HABalt
$$
Indeed, denote by $H'_{A, B}$ the right-hand side of \ref e.HABalt/.
Suppose that $u \in H_{A, B}$ and
write $u = u_1+u_2+u_3$ with $u_1 \in H \cap A^\perp$, $u_2 \in [A]$, and $u_3
\in [B]$. Since $u_1 \in H$, it follows that $u_1+u_3 \in H+[B]$. Since $u_1
\in A^\perp$ and $u_3 \in [B] \subseteq A^\perp$, we obtain $u_1+u_3 \in
(H+[B]) \cap A^\perp$. Since $u_2 \in [A] \subseteq B^\perp$ and $u \in
B^\perp$, we have that $u_1+u_3 = u-u_2 \in B^\perp$, whence $u_1+u_3 \in
\big( (H+[B]) \cap A^\perp \big) \cap B^\perp = ( H + [B] ) \cap (A \cup
B)^\perp$. Since $u_2 \in [A]$, it follows that $u \in H'_{A, B}$.

Conversely, suppose that $u \in H'_{A, B}$ and write $u = u_1+u_2+u_3$ with
$u_1 \in H$, $u_2 \in [B]$, $u_3 \in [A]$, and $u_1+u_2 \in [A \cup B]^\perp$.
Since $u_1+u_2 \in A^\perp$ and $u_2 \in A^\perp$, it follows that $u_1 \in
A^\perp$, whence $u_1 \in H \cap A^\perp$ and so
$u \in (H \cap A^\perp) + [A \cup B]$. Since $u_1+u_2 \in
B^\perp$ and $u_3 \in B^\perp$, we also have $u \in B^\perp$, whence $u \in
H_{A, B}$, as desired.

\procl c.commute
For any closed subspace $H$ and $A, B, C, D \subseteq E$ with $(A \cup B) \cap
(C \cup D) = \emptyset$, we have 
$$ 
(H_{A, B})_{C, D} = H_{A \cup C, B \cup D}
%(H_{A, B})_{C, D} = (H_{C, D})_{A, B}
\,.
$$
\endprocl

\proof
In showing \ref e.HABalt/,
we have shown that $H_{A, B} = (H_{A, \emptyset})_{ \emptyset, B} = (H_{
\emptyset, B})_{A, \emptyset}$.
Also, we have 
\begineqalno
(H_{A, \emptyset})_{B, \emptyset}
&=
\Big(\big((H \cap A^\perp) +[A]\big)\cap B^\perp \Big) +[B]
\cr&=
\Big(\big((H \cap A^\perp) +[A \setminus B]\big)\cap (B \setminus A)^\perp
\Big) +[B]
\cr&=
\big((H \cap (A \cup B)^\perp) +[A \setminus B]\big) +[B]
\cr&=
(H \cap (A \cup B)^\perp) +[A \cup B]
=
H_{A \cup B, \emptyset}
\,,
\endeqalno
whence also $(H_{ \emptyset, A})_{ \emptyset, B} = H_{ \emptyset, A \cup B}$
by duality.
Therefore, 
\begineqalno
(H_{A, B})_{C, D} 
&= (((H_{A, \emptyset})_{ \emptyset, B})_{
\emptyset, D})_{C, \emptyset}
= ((H_{A, \emptyset})_{ \emptyset, B \cup D})_{C, \emptyset}
\cr&= ((H_{A, \emptyset})_{C, \emptyset})_{ \emptyset, B \cup D}
= (H_{A \cup C})_{ \emptyset, B \cup D}
= H_{A \cup C, B \cup D}
\,.
\Qed
\cr
\endeqalno

We remark that when $A \cup B$ is finite, the subspace $H_{A, B}$ is closed,
since the sum of two closed subspaces, one of finite dimension, is always
closed (\ref b.Halmos:prob/, Problem 13) and since the intersection of closed
subspaces is always closed.

It follows by induction from \ref p.minors/ that for {\it finite} $E$,
if $\P^H[A \subseteq \ba, B
\cap \ba = \emptyset] > 0$, then 
$$
\P^H(\;\cbuldot \mid A \subseteq \ba, B \cap \ba = \emptyset)
=
\P^{H_{A, B}}(\;\cbuldot\;)
\,.
\label e.HABcond
$$

For a set $K
\subseteq E$, let $\F(K)$ denote the $\sigma$-field 
generated by the events $ \{ e \in \ba \} $ for $e \in K$.
We shall say that the events in $\F(K)$ are {\bf measurable with respect to}
$K$ and likewise for functions that are measurable with respect to $\F(K)$.
We also say that an event or a function 
that is measurable with respect to $E \setminus
\{ e \}$ {\bf ignores} $e$.
%$$
%A \cup \{e\} \in \ev A \iff A \in \ev A
%\,.
%$$
Thus, an event $\ev A$ ignores $e$ iff for all $\u \in S_{\A}$, also
$\u \wedge e \in S_{\ev A}$ and $\u \vee e \in S_{\ev A}$.
We say that a probability measure $\P$ on $2^E$ has {\bf negative
associations} if
$$
\P(\A_1 \cap \A_2) \le \P(\A_1) \P(\A_2)
$$
whenever $\A_1$ and $\A_2$ are increasing events that are measurable with
respect to 
complementary subsets of $E$.
A function $f : 2^E \to \R$ is called {\bf increasing} if for all $A \in
2^E$ and all $e \in E$, we have $f(A \cup \{ e \}) \ge f(A)$.
%We say that a function $f$ {\bf ignores} $e$ if for all $A \in 2^E$, we
%have $f(A \cup \{ e \}) = f(A)$.
It is well known and not hard to see that $\P$ has negative associations iff
for any pair $f_1$, $f_2$ of increasing functions that are measurable with
respect to complementary subsets of $E$,
$$
\E[f_1 f_2] \le \E[f_1] \E[f_2] 
\,.
\label e.negass
$$
In this case, for any collection $f_1, f_2, \ldots, f_n$ of increasing
{\it nonnegative\/}
functions that are measurable with respect to pairwise disjoint
subsets of $E$, we have
$$
\E[f_1 f_2 \cdots f_n] \le \E[f_1] \E[f_2] \cdots \E[f_n]
\,.
\label e.negass-multi
$$
This is shown by an easy induction argument.
One could replace ``increasing" by ``decreasing" just as well.

The following result was proved by \ref b.FedMih/, Lemma 3.2, for
the uniform measure on bases of ``balanced" matroids, a class of matroids
including regular ones.
%{An example of a matroid that is not balanced is given on p.~495 of
%Seymour and Welsh (1975). This is referred to in Feder and Mihail.}
%{An example of balanced nonregular matroid is $U_{2, 4}$, the
%uniform matroid of all 2-element subsets of a 4-element ground set. This is
%also a real matroid: just take 4 vectors in $\R^2$ in general position.}

\procl t.FM 
Let $E$ be finite.
If $\A$ is an increasing event that ignores $e$ and $\P^H[e \in \ba] > 0$,
then $\P^H(\A \mid e \in \ba) \le \P^H(\A)$. More generally, $\P^H$ has
negative associations.
\endprocl

\proof
We have $\P^H(\A \mid e\in \ba) = \P^{(H \cap e^\perp)+[e]}(\A) = \P^{(H \cap
e^\perp)}(\A)$ since the effect of having $e$ in a subspace is simply to
make $e \in \ba$ with probability 1 and since $\A$ ignores $e$. Since $H
\cap e^\perp \subseteq H$, it follows from \ref t.dominate/ that $\P^{(H \cap
e^\perp)}(\A) \le \P^H(\A)$, which proves the first assertion.

For the more general statement, we follow the method of proof of \ref
b.FedMih/.
This uses induction on the cardinality of $E$.
The case $|E| = 1$ is trivial.
Let $r := \dim H$.
Given $\A_1$ and $\A_2$ as specified with $\P^H(\A_1), \P^H(\A_2) > 0$,
%there is an
%element $e$ on which $\A_1$ depends (and so $\A_2$ ignores) such that
%$\A_1$ is positively correlated with $e$, since $\A_1$ is negatively
%correlated with those elements that $\A_1$ ignores and since the number of
%elements in $\ba$ is constant: If 
we have
$$
\sum_{e \in E} \P^H[e \in \ba] = \E^H\big[|\ba|\big] = r
= \E^H\big[|\ba| \bigm| \A_1\big] = \sum_{e \in E} \P^H[e \in \ba \mid \A_1]
$$
since $|\ba| = r$ $\P^H$-a.s.
By the preceding paragraph,
$\P^H[e \in \ba \mid \A_1] \le \P^H[e \in \ba]$ for those $e$ ignored
by $\A_1$. It follows that the opposite inequality holds for some $e$
not ignored by $\A_1$ and with $\P^H[e \in \ba] > 0$.
Fix such an $e \in E$.
Note that $\A_2$ ignores $e$.
We have $\P^H(\A_1 \mid e \in \ba) \ge \P^H(\A_1 \mid e \notin \ba)$, where
the right-hand side is defined to be 0 if $\P[e \notin \ba] = 0$.
Now
$$
\P^H(\A_1 \mid \A_2)
=
\P^H[e \in \ba \mid \A_2] \P^H(\A_1 \mid \A_2, e \in \ba)
+
\P^H[e \notin \ba \mid \A_2] \P^H(\A_1 \mid \A_2, e \notin \ba)
\,.
\label e.1
$$
The induction hypothesis and \ref p.minors/ imply
that \ref e.1/ is at most
$$
\P^H[e \in \ba \mid \A_2] \P^H(\A_1 \mid e \in \ba)
+
\P^H[e \notin \ba \mid \A_2] \P^H(\A_1 \mid e \notin \ba)
\,;
\label e.2
$$
this is because by \ref p.minors/,
the two measures conditioned on whether or not $e$ lies in
$\ba$ can each be regarded as measures arising from orthogonal
projections on $\ell^2(E
\setminus \{ e \})$ and because $\A_1$ and $\A_2$ each transform to
increasing events in $2^{E \setminus \{ e \}}$.
By what we have proved in the first paragraph, we have that $\P^H[e \in \ba
\mid \A_2] \le \P^H[e \in \ba]$ and we have chosen $e$ so that
$\P^H(\A_1 \mid e \in \ba) \ge \P^H(\A_1 \mid e \notin \ba)$.
Therefore, \ref e.2/ is at most
$$
\P^H[e \in \ba] \P^H(\A_1 \mid e \in \ba)
+
\P^H[e \notin \ba] \P^H(\A_1 \mid e \notin \ba)
=
\P^H(\A_1)
\label e.3
$$
since the quantity in \ref e.3/ minus that in \ref e.2/ equals 
$$
\big(\P^H[e \in \ba] - \P^H[e \in \ba \mid \A_2]\big)
\big(\P^H(\A_1 \mid e \in \ba) - \P^H(\A_1 \mid e \notin \ba)\big)
\ge 0
\,.
$$
This completes the induction step and hence the proof.
\Qed

It seems that there should be an elegant proof of \ref t.FM/ along the lines
of the proof of \ref t.dominate/, but we were not able to find one.
However, the special case of \ref t.FM/ saying that 
$$
\P^H[A \cup B \subset \ba] \le \P^H[A \subset \ba] \P^H[B
\subset \ba]
$$
is an immediate consequence of \ref e.whitney/ for $\u := \bigwedge_{e \in A}
P_H e$ and $\v := \bigwedge_{e \in B} P_H e$.
This special case was also proved in Proposition 2.7 of \ref b.ShiTak:II/.

Combining \ref t.FM/ with Equation \ref e.HABcond/, we obtain the following
result:

\procl c.condnegass
If $E$ is finite, $A$ and $B$ are disjoint subsets of $E$,
and $H$ is a subspace of $\H$, then the conditional probability measure
$\P^H(\;\cbuldot \mid A \subseteq \ba, B \cap \ba = \emptyset)$ has negative
associations.
\endprocl

In the language of \ref b.Pemantle:negass/, this corollary says that the
measure $\P^H$ is conditionally negatively associated (CNA). 
In fact, it enjoys the
strongest property studied by \ref b.Pemantle:negass/, namely, CNA with
so-called external fields. This follows by \ref r.weight/.

\ref b.Soshnikov:fluct/ proves a very general central limit theorem for
determinantal probability measures (and determinantal point processes).
Many theorems for independent random variables are known to hold for
negatively associated random variables, especially for a stationary
collection of negatively random variables indexed by $\Z^d$. For example,
see \ref b.Newman:asympt/, \ref b.ShaoSu:LIL/, \ref b.Shao:comparison/,
\ref b.ZhangWen:weak/, \ref b.Zhang:Strassen/, and the references therein.
We shall merely state one of the easiest of these theorems, as it will
prove useful to us later:
%It is known that
%the usual proof of sub-Gaussian bounds for bounded independent random
%variables applies for negatively associated bounded random variables.
%Therefore, we obtain the following bounds:

\procl c.subG
Let $E$ be finite, $H$ be a subspace of $\H$, and $A \subseteq
E$. Write $\mu := \E^H[|\ba \cap A|]$. Then for any $a > 0$, we have 
$$
\P^H\bigg[\Big| |\ba \cap A| - \mu \Big|  \ge a\bigg] \le 2 e^{-2 a^2/|A|}
\,.
\label e.subG
$$
\endprocl

%\proof
Essentially the same observation is made by \ref b.DubRan/.
%For fixed real $\lambda$, the functions on $2^E$
%$$
%B \mapsto \exp \Big\{ \lambda \II{e \in B} \Big\}
%%B \mapsto \exp \Big\{ \lambda \big( \II{e \in B} - \P^H[e \in \ba]\big) \Big\}
%\label e.fne
%$$
%for $e \in A$ are simultaneously all increasing or all decreasing nonnegative
%functions.
%They also ignore complementary subsets, as the function \ref e.fne/
%ignores all $e' \ne e$.
%Hence the standard proof, with the addition of \ref e.negass-multi/ to replace
The standard proof, with the addition of \ref e.negass-multi/ to replace
independence, applies.
See, e.g., \ref b.AS:book/, Corollary A.1.7, for such a proof for
independent random variables.
%\Qed

Recall Kirchhoff's theorem that the current $Y(e,e)$ is the probability that
$e$ lies in a uniform spanning tree of a graph $G$.
Kirchhoff extended this equation to express the entire current vector in a
network by a sum over spanning trees; see \ref b.Thomassen:res/ for a
statement and short combinatorial proof of this result.
Such a combinatorial interpretation and proof extends to the uniform
measure on bases of regular matroids, but not, as far as we can tell, to
our general case.
We do, however, have the following extension,
in which \ref e.expzeta/
is related to, but seems not to follow from, \ref b.Maurer/, Theorem 1.

\procl p.kirchproj
Let $E$ be finite.
For any $v \in \H$ and $B \in \B$, define $a_v(e, B) = a_v(e, B; H)$ by
$$
P_H v = \sum_{e\in B} a_v(e, B) P_H e
\,;
\label e.defa_v
$$
this uniquely determines the coefficients $a_v(e, B)$ since the vectors
$\Seq{P_H e \st e \in B}$ form a basis of $H$. Put
$$
\zeta^v_B := \zeta^v_B(H) := \sum_{e\in B} a_v(e, B) e
\,.
$$
Then 
$$
P_H v = \E^H \zeta^v_\ba
\,.
\label e.expzeta
$$
In addition, if $S \subseteq F \subset E$ and $e \notin F$, then
$$
P_\HFS e = P_{[F]}^\perp \E^H[\zeta^e_\ba(H) \mid \ba \cap F = S]
\,,
\label e.expzetaHFS
$$
where
$$
\HFS := \big(H + [F \setminu S] \big) \cap [F]^\perp
\,.
\label e.HFSdef
$$
\endprocl

\comment{Does this have an explicit form for group representations from the
canonical decomposition?}

In order to understand this proposition better, we interpret it for the
graphical matroid of a finite graph $G$ and $v = e \in\edges$.
In this case, we have $a_e(f, B) \in \{ 0, 1 \}$ for every edge $f$ and
$\zeta^e_B$ is the path in $B$ from the tail of $e$ to the head of $e$.
Equation \ref e.expzeta/ is then the theorem of \ref b.Kirchhoff/ alluded to
above.
Equation \ref e.expzetaHFS/ is the same thing applied to a minor of $G$.

\proof 
%Since $P_H \E^H \zeta^v_\ba = \E^H P_H \zeta^v_\ba = \E^H P_H v = P_H v$, it
%suffices for \ref e.expzeta/ to show that $\E^H \zeta^v_\ba \in H$.
%
For $B \in \B$, let $\theta_B = \bigwedge_{i=1}^r e_i$.
For any $j$, if we take the wedge product of \ref e.defa_v/ with
$\bigwedge_{i < j} P_H e_i$
on the left and $\bigwedge_{i > j} P_H e_i$ on the right and use \ref
l.projection/, we obtain
$$
P_{\Ext(H)} \left(
\bigwedge_{i < j} e_i \wedge v \wedge \bigwedge_{i > j} e_i \right)
=
a_v(e_j, B) P_{\Ext(H)} \theta_B
\,,
$$
which is the same as 
$$
P_{\Ext(H)} [(\theta_B \vee e_j) \wedge v] 
=
a_v(e_j, B) P_{\Ext(H)} \theta_B
\,.
$$
Therefore,
$$
a_v(e_j, B) = \bigip{(\theta_B \vee e_j) \wedge v, \xi_H}/\ip{\theta_B, \xi_H}
\,.
\label e.a_v=
$$
Since $P^H[B] = \ip{\xi_H, \theta_B}\ip{\theta_B, \xi_H}$
by \ref e.xiHpr/ and $\theta_B \vee
e = 0$ for $e \notin B$, it follows that
\begineqalno
\E^H \zeta^v_\ba
&=
\sum_B \ip{\xi_H, \theta_B} \sum_{e \in B} \bigip{(\theta_B \vee e) \wedge
v, \xi_H} e
=
\sum_B \ip{\xi_H, \theta_B} \sum_{e \in E} \bigip{\theta_B, (\xi_H \vee
v) \wedge e} e
\cr&=
\sum_{e \in E} e \sum_B \ip{\xi_H, \theta_B} \bigip{\theta_B, (\xi_H \vee
v) \wedge e} 
=
\sum_{e \in E} e \bigip{\xi_H, (\xi_H \vee v) \wedge e} 
\cr&\qquad\qquad\hbox{[since $\theta_B$ are orthonormal and $\xi_H$ lies in
their span]}
\cr&=
\sum_{e \in E} \bigip{\xi_H \vee e, \xi_H \vee v} e 
=
\sum_{e \in E} \ip{P_H v, e} e
=
P_H v
\endeqalno
by \ref e.reverse/.
%This gives that for any $w \in \H$,
%$$
%\bigip{\E^H \zeta^v_\ba, w}
%=
%\ip{\xi_H \vee w, \xi_H \vee v}
%\,.
%$$
%If $w \perp H$, then $\xi_H \vee w = 0$ (complete $w$ to an orthonormal
%basis of $\H$ containing a basis of $H$ to see this), whence $\bigip{\E^H
%\zeta^v_\ba, w} = 0$. That is, $\E^H \zeta^v_\ba \in H$.
This proves \ref e.expzeta/.

To show \ref e.expzetaHFS/, note that
$$
H_{S, F \setminus S} = 
\HFS + [S]
\,.
\label e.defH^F_S
$$
By \ref e.HABcond/, we have that $\P^H(\; \cbuldot \mid \ba
\cap F = S) = \P^{\HSFS}(\; \cbuldot \;)$.  In fact, induction from \ref
p.minors/ shows that
$\ip{\theta_B, \xi_{\HSFS}}/\ip{\theta_B, \xi_H}$ does not depend on $B$ so
long as $B \cap F = S$; this quotient is 0 if $B \cap F \ne S$. 
Now $B \cap F = S$ and $e, e' \notin F$ imply that $\big( (B \setminu e')
\cup e \big) \cap F = S$.
Inspection of \ref e.a_v=/ shows, therefore, that if $B \cap F = S$ and $e
\notin F$, then $a_e(e', B; H) = a_e(e', B; {\HSFS})$ for $e' \in B\setminu
F$, and thus $P_{[F]}^\perp \zeta^e_B(H) = \sum_{e' \in B \setminu F}
a_e(e', B; H) e' = \zeta^e_B(\HSFS)$.
Therefore,
\begineqalno
P_{[F]}^\perp \E^H[ \zeta^e_B(H) \mid B \cap F = S]
&=
P_{[F]}^\perp \E^H[ \zeta^e_B(\HSFS) \mid B \cap F = S]
\cr&=
P_{[F]}^\perp \E^{\HSFS}[ \zeta^e_B(\HSFS) ]
=
P_{[F]}^\perp P_{\HSFS} e
\cr&=
P_{\HFS} e
\,.
\Qed
\endeqalno

\bsection{The Infinite Case}{s.infinite}

Most interesting results about uniform spanning trees and forests arise in
the setting of infinite graphs.
In this section and the following one, we shall see what can be
accomplished for general determinantal probability measures on
infinite ground sets $E$.
We still restrict ourselves in this section to orthogonal projections. In the
next section, we extend both the finite and infinite cases to positive
contractions.
%We first show how to construct determinantal probability measures
%on infinite sets and give some of their properties.

Let $|E| = \infty$ and consider first a finite-dimensional subspace $H$ of
$\H$. Let us order $E$ as $\{e_i \st i \ge 1\}$. Define $H_k$ as the image of
the orthogonal projection of $H$ onto the span of $\{e_i \st 1 \le i \le k\}$.
By considering a basis of $H$, we see that $P_{H_k} \to P_H$ in the strong
operator topology (SOT), i.e., for all $v \in \H$, we have
$\|P_{H_k} v - P_H v\| \to 0$ as $k \to\infty$. 
We shall write $H_k \SOTto H$ for $P_{H_k} \to P_H$ in the SOT.
It is also easy to
see that if $r := \dim H$, then $\dim H_k = r$ for all large $k$ and, in
fact, $\xi_{H_k} \to \xi_H$ in the usual norm topology. It follows that
\ref e.genprs/ holds for this subspace $H$ and for any finite $A_1,
A_2 \subset E$.

Now let $H$ be an infinite-dimensional (closed) subspace of $\H$. 
It is well known that if $H_n$ are (closed) subspaces of $\H$ with $H_n
\uparrow H$ (meaning that $H_n \subseteq H_{n+1}$ and $\bigcup H_n$ is dense
in $H$) or with $H_n \downarrow H$ (meaning that $H_n \supseteq H_{n+1}$ and
$\bigcap H_n = H$), then $H_n \SOTto H$.
[The proof follows
immediately from writing $H$
as the orthogonal direct sum of its subspaces $H_{n+1} \cap H_n^\perp$ in the
first case and then by duality in the second.]
Choose an increasing sequence of finite-dimensional subspaces $H_k \uparrow
H$.
Since $H_k \SOTto H$, we have 
$$
\hbox{for all finite sets } A \quad
\det (P_{H_k} \restrict A) \to \det (P_H \restrict A)
\,,
\label e.detgenprs
$$
whence $\P^{H_k}$ has a weak${}^*$ limit that we
denote $\P^H$ and that satisfies
\ref e.genprs/.
{}From this construction, we see that
the statement of \ref t.dominate/ is valid for two possibly
infinite-dimensional subspaces, one contained in the other:

\procl t.dominate-infinite
Let $E$ be finite or infinite and let
$H_1 \subset H_2$ be closed subspaces of $\H$. Then $\P^{H_1} \preccurlyeq
\P^{H_2}$.
\endprocl

We also note that for {\it any\/} sequence of subspaces $H_k$, if $H_k \SOTto
H$, then $\P^{H_k} \to \P^H$ weak${}^*$ because \ref e.detgenprs/ then holds.

%Given any finite subset $F
%\subset E$, the probabilities $\P^{H_k}[F \subset \ba]$ are increasing in
%$k$ by \ref t.dominate/. Thus, they have a limit, $\P^H[F \subset \ba]$. These
%probabilities determine the probabilities of all cylinder events, hence
%define the probability measure $\P^H$. Most properties that hold for the
%finite-dimensional case extend in a similar manner. 

Establishing a conjecture of \BLPSusf,
\ref b.Morris/ proved that on any network $(G, w)$ (where $G$ is the underlying
graph and $w$ is the function assigning conductances, or weights, to the
edges), for $\wsf(G, w)$-a.e.\ forest $\fo$ and for every component tree $T$
of $\fo$, the $\wsf$ of $(T, w \restrict T)$ equals $T$ a.s.
This suggests the following extension. Given a
subspace $H$ of $\H$ and a set $B \subseteq E$, the subspace of $[B]$
``most like" or ``closest to" $H$ is the closure of the image of $H$ under the
orthogonal projection $P_{[B]}$; we denote this
subspace by $H_B$. For example, if $H = \STAR(G)$,
then $H_B = \STAR(B)$ since
for each $x \in V(G)$, we have $P_{[B]} (\star_x^G) = \star_x^B$.
To say that $\P^{H_B}$ is concentrated on $ \{ B \}$ is the same as to say
that $H_B = [B]$. This motivates the following theorem and shows how it is
an extension of Morris's theorem.

\procl t.morris
For any closed subspace $H$ of $\H$, we have $H_{\ba} = [\ba]$ $\P^H$-a.s.
\endprocl

After the proof of this theorem, we shall give some applications.
In order to establish \ref t.morris/, we shall use several short lemmas.
The overall proof is modelled on the original (second)
proof by Morris.\ftnote{*}
{His proof is much shorter than our proof
as he could rely on known facts about
electrical networks. In particular, the relationship of our proof to that
of Morris may not be so apparent. The part of our proof that is easiest to
recognize is \ref l.rayleigh/, which is a version of Rayleigh's
monotonicity principle for electric networks.}

\procl l.HB
For any closed subspace $H$ and any $B \subseteq E$, we have $H_B
= \overline{(H + B^\perp) \cap [B]} = \overline{H + B^\perp} \cap [B]$.
\endprocl

\proof
Define $H'_B := (H + B^\perp) \cap [B]$
and $H''_B := \overline{H + B^\perp} \cap [B]$.
We shall show that $H_B \subseteq \overline{H'_B} \subseteq H''_B \subseteq
H_B$.
Given $u \in H$, write $u = u_1+u_2$ with $u_1 \in [B]$ and $u_2 \in B^\perp$.
Then $u_1 = u - u_2 \in H + B^\perp$, and so $P_{[B]} u = u_1 \in H'_B$.
Therefore, $H_B \subseteq \overline{H'_B}$.
It is clear that $H'_B \subseteq H''_B$ and that $H''_B$ is closed, whence
$\overline{H'_B} \subseteq H''_B$.
Finally, given $u \in H''_B$, write $u = \lim_{n
\to\infty} u^{(n)}$, with $u^{(n)} = u^{(n)}_1+u^{(n)}_2$, where $u^{(n)}_1
\in H$ and $u^{(n)}_2 \in B^\perp$.
Since $u \in [B]$, we have 
$$
u = P_{[B]} u = P_{[B]} \lim_{n \to\infty} u^{(n)}
= \lim_{n \to\infty} P_{[B]} u^{(n)} = \lim_{n \to\infty} P_{[B]} u^{(n)}_1 \in H_B
\,.
\Qed
$$

\procl l.HBplus
For any closed subspace $H$ and any $B \subseteq E$, we have
$H_B + B^\perp = \overline{H + B^\perp}$.
\endprocl

\proof
For any $u \in H$, write $u = u_1+u_2$ with $u_1 \in [B]$ and $u_2 \in
B^\perp$.
We have that $u + B^\perp = u_1 + B^\perp = P_{[B]} u + B^\perp$.
Since the closure of $ \{ u + B^\perp \st u \in H \}$ equals $\overline{H +
B^\perp}$ and the closure of
$\{ P_{[B]} u + B^\perp \st u \in H \}$ equals $H_B + B^\perp$, the result
follows.
\Qed

\procl l.rayleigh
For any closed subspace $H$ and
$B \subseteq E$, if $v \in [B]$, then $\|P_{H_B} v\| \ge \|P_H v\|$.
\endprocl

\proof
We have $P_{H_B} v = P_{H_B + B^\perp} v$ since $v \in [B]$.
By \ref l.HBplus/, it follows that $P_{H_B} v = P_{\overline{H + B^\perp}} v$.
Since $H \subseteq \overline{H + B^\perp}$, it follows that $P_H v = P_H
P_{H_B} v$, whence $\|P_H v\| \le \|P_{H_B} v\|$, as desired.
\Qed

\procl r.rayleigh
More generally,
if $v \in [B]$ and $B \subseteq C \subseteq E$, then $\|P_{H_B} v\| \ge
\|P_{H_C} v\|$.
Indeed,
since $P_{[B]} H = P_{[B]} P_{[C]} H$, we have $H_B = (H_C)_B$.
Thus, the inequality follows from applying \ref l.rayleigh/ to the space
$H_C$, rather than to $\H$.
\endprocl

\procl l.goingto0
Fix $e \in E$.
Let $B \subseteq E$ and $H$ be a closed subspace of $[B]$.
Suppose that $F_n \subseteq E$ form an increasing sequence of
sets with union $E \setminus \{ e \}$.
If $e \notin H$, then 
$\lim_{n \to\infty} \|P_{H_{F_n \cap B, F_n \setminus B}} e\| = 0$.
\endprocl

\proof
Since $H \subseteq [B]$, we have by \ref e.HAB/ that
$$
H_{F_n \cap B, F_n \setminus B}
=
(H \cap (F_n \cap B)^\perp) + [F_n \cap B]
\,,
$$
whence 
$$
P_{H_{F_n \cap B, F_n \setminus B}} e
=
P_{H \cap (F_n \cap B)^\perp} e + P_{[F_n \cap B]} e
=
P_{H \cap (F_n \cap B)^\perp} e 
$$
since $e \notin F_n$.
Also, 
\begineqalno
\bigcap (H \cap (F_n \cap B)^\perp) 
&=
H \cap (\bigcup F_n \cap B)^\perp
=
H \cap (B \setminus \{ e \})^\perp
\cr&=
H \cap ([e] + [E \setminus B])
=
H \cap [e]
\cr
\endeqalno
since $H \subseteq
[B]$, whence $H \cap (F_n \cap B)^\perp \SOTto H \cap [e] = 0$
if $e \notin H$.
\Qed

\proofof t.morris
Fix $e \in E$.
Let $\A_1 := \{B \in 2^E \st e \in B \}$ and $\A_2 := \{B \in 2^E \st e \notin
H_B \}$.
We want to show that $\P^H(\A_1 \cap \A_2) = 0$.
Let $E \setminus \{ e \} = \bigcup F_n$ for increasing finite sets $F_n$.
By \ref e.HABcond/, we have 
$$
\P^H\big(\A_1 \mid \F(F_n)\big)(B) = 
\P^{H_{F_n \cap B, F_n \setminus B}}(\A_1)
=
\|P_{H_{F_n \cap B, F_n \setminus B}} e\|^2
\,.
$$
Now 
$$
(H_{F_n \cap B, F_n \setminus B})_B
=
\overline{(H_{F_n \cap B, F_n \setminus B})_{ \emptyset, E \setminus B}}
\label e.1
$$
by \ref l.HB/.
Write $H_B^n := (H_B)_{F_n \cap B, F_n \setminus B}$.
Combining \ref e.1/ with \ref c.commute/ and \ref l.HB/, we obtain 
$$
(H_{F_n \cap B, F_n \setminus B})_B
=
\overline{(H_{ \emptyset, E \setminus B})_{F_n \cap B, F_n \setminus B}}
\subseteq
\overline{(H_B)_{F_n \cap B, F_n \setminus B}}
=
H_B^n
$$
since $F_n$ is finite.
If $B \in \A_1$, then we may apply 
\ref l.rayleigh/ to $H_{F_n \cap B, F_n \setminus B}$ and obtain
$$
\P^H\big(\A_1 \mid \F(F_n)\big)(B)
=
\|P_{H_{F_n \cap B, F_n \setminus B}} e\|^2
\le
\|P_{(H_{F_n \cap B, F_n \setminus B})_B} e\|^2
\le
\|P_{H_B^n} e\|^2
\,,
$$
and this tends to 0
as $n \to\infty$ if $B \in \A_2$ by \ref l.goingto0/ (applied to $H_B
\subseteq [B]$).
On the other hand, $\P^H\big(\A_1 \mid \F(F_n)\big) \to \P^H\big(\A_1
\mid \F(E \setminus \{ e \})\big)$ a.s.\ by L\'evy's martingale convergence
theorem. Therefore, $\P^H\big(\A_1 \mid \F(E \setminus \{ e \})\big) = 0$ a.s.\
on the event $\A_1 \cap \A_2$. 
Write $\A_3 := \Big\{ B \st \P^H\big(\A_1 \mid \F(E \setminus \{ e \})\big)(B)
= 0 \Big\}$, an event that lies in $\F(E \setminus \{ e \})$ and contains
$\A_1 \cap \A_2$.  We have
\begineqalno
\P^H(\A_1 \cap \A_2) 
&\le \P^H( \A_1 \cap \A_3)
=
\E^H\Big[\P^H\big( \A_1 \cap \A_3 \mid \F(E \setminus \{ e \})\big)\Big]
\cr&=
\E^H\Big[\P^H\big( \A_1 \mid \F(E \setminus \{ e \})\big) \I{\A_3} \Big]
=
0
\,.
\endeqalno
Since this holds for each $e \in E$ and $E$ is countable, the
theorem follows.
\Qed

We now give some applications of \ref t.morris/.
For the special case of $H := \CYCLE(G)^\perp$, we get the following result.

\procl c.FSFcor 
Let $(G, w)$ be a network. For the free spanning forest $\fo$, we have
a.s.\ $(\CYCLE(G)^\perp)_\fo = [\fo]$.
\endprocl

% This is not true for all forests in $G$; e.g., $G := T \times \Z$ and $\fo
% =$ union of $T$ fibers.

This is a nontrivial theorem
about the free spanning forest, different from the trivial statement that
the free spanning forest of the free spanning forest is itself, since
$(\CYCLE(G)^\perp)_\fo \subsetneq \CYCLE(\fo)^\perp$ except in degenerate
cases.
The dual statement is easier to interpret:

\procl c.dualFSF
Let $(G, w)$ be a network. For the free spanning forest $\fo$, we have
that a.s.\ the free spanning forest of the contracted graph $G/\fo$ is
concentrated on the empty set.
\endprocl

Note that $G/\fo$ may have vertices of infinite degree.

\proof
By duality, the free spanning forest has the law of $\edges \setminus \ba$
when $\ba$ has the law of $\P^{\CYCLE(G)}$.
We can naturally identify $(\CYCLE(G))_\ba$ with $\CYCLE(G/\ba)$.
Thus, the fact that $(\CYCLE(G))_\ba = [\ba]$ $\P^{\CYCLE(G)}$-a.s.\ gives the
statement of the corollary by duality again.
\Qed

We next give a dual form of \ref t.morris/, a form that is a very
natural property for infinite matroids with respect to a probability
measure $\P^H$. Then we give some applications of this reformulation.
The duality that we now use is the following. 

\procl l.perpperp
Let $H_1$ and $H_2$ be any two closed subspaces of a Hilbert space with
corresponding orthogonal projections $P_1, P_2$ and coprojections $P_1^\perp,
P_2^\perp$.
Then $\overline{P_1 H_2} = H_1$ iff $\overline{P_2^\perp H_1^\perp} =
H_2^\perp$.
\endprocl

\proof
By symmetry, it suffices to show that if $\overline{P_1 H_2} \ne H_1$, then
$\overline{P_2^\perp H_1^\perp} \ne H_2^\perp$.
Now if $\overline{P_1 H_2} \ne H_1$, then there exists a non-0 vector
$u \in H_1 \cap (P_1 H_2)^\perp$.
Fix such a $u$.
For all $v \in H_2$, we have $0 = \ip{u, P_1 v} = \ip{P_1 u, v} = \ip{u,
v}$, so that $u \in H_2^\perp$.
Therefore, for any $w \in H_1^\perp$, we have $\ip{u, P_2^\perp w}
= \ip{P_2^\perp u, w} = \ip{u, w} = 0$, so that $u \perp
P_2^\perp(H_1^\perp)$.
Hence, $u \in H_2^\perp \cap \left(P_2^\perp(H_1^\perp)
\right)^\perp$, so that $\overline{P_2^\perp H_1^\perp} \ne H_2^\perp$.
\Qed

Recall that when $H$ is finite dimensional, $\P^H$ is supported by those
subsets $B \subseteq E$ that project to a basis of $H$ under $P_H$.
(Strictly speaking, we have shown this only when $E$ is finite. However, the 
definition shows that it is true when $E$ is infinite as well, provided $H$
is still finite-dimensional.)
The following theorem extends this to the infinite setting insofar
as a basis is a spanning set. (The other half of being a basis,
minimality, does not hold in general, even for
the wired spanning forest of a tree, as shown by the examples in \ref
b.HL:change/.)
% Perhaps invariance would make it hold. Then we would get bilaterally
% deterministic Bernoulli shifts. But I think not.

\procl t.basis
For any closed subspace $H \subseteq \H$, we have $\overline{[P_H \ba]} = H$
$\P^H$-a.s.
\endprocl

% By separability, it is enough to show a single vector in $H$ is in the
% closure of $[P_H \ba]$ a.s.

\proof
%For the moment, fix $B \subseteq E$.
%Suppose that $u \in H \cap [P_H B]^\perp$.
%Then for all $e \in B$, we have $0 = \ip{u, P_H e} = \ip{P_H u, e} = \ip{u,
%e}$, so that $u \in [E \setminus B]$.
%Therefore, for any $v \in H^\perp$, we have $\ip{u, P_{[E \setminus B]} v}
%= \ip{P_{[E \setminus B]} u, v} = \ip{u, v} = 0$, so that $u \perp
%(H^\perp)_{[E \setminus B]}$.
%Hence, $u \in [E \setminus B] \cap \left((H^\perp)_{[E \setminus
%B]}\right)^\perp$, so that if $(H^\perp)_{[E \setminus B]} = [E \setminus
%B]$, then $u = 0$.
%
%Now by \ref t.morris/ applied to $H^\perp$ and rewritten by using
%\ref c.dualrep/, we have $(H^\perp)_{[E \setminus \ba]} = [E \setminus
%\ba]$ $\P^H$-a.s.
%Thus, our calculations of the preceding paragraph show that $H \cap [P_H
%\ba]^\perp = 0$ $\P^H$-a.s., which is equivalent to the desired result.
According to \ref l.perpperp/, $\overline{[P_H B]} = H$ is equivalent to 
$(H^\perp)_{[E \setminus B]} = [E \setminus B]$.
If we apply \ref t.morris/ to $H^\perp$ and rewrite the conclusion by using
\ref e.dualrep/, we see that this holds for $\P^H$-a.e.\ $B$.
\Qed

\procl r.betterproof
This reasoning shows
that \ref t.morris/ could be deduced from an alternative proof of \ref
t.basis/. Thus, it would be especially worthwhile to find a simple direct
proof of \ref t.basis/.
\endprocl

Our first application of \ref t.basis/ is for $E = \Z$. 
% but it has an analogous statement for any countable abelian group.
Let $\T := \R/\Z$ be the unit circle equipped with unit Lebesgue measure. For
a measurable function $f : \T \to \C$ and an integer $n$, the {\bf Fourier
coefficient} of $f$ at $n$ is 
$$
\widehat f(n) := \int_\T f(t) e^{-2\pi i n t} \,dt
\,.
$$
%Given a Lebesgue measurable set $A \subset \T$ and a set $S \subseteq \Z$, say
%that $S$ is a {\bf restriction set for} $A$ if the set $\{\widehat f
%\restrict S \st f \in L^2(A)\}$ is dense in $\ell^2(S)$.
%(Here, $L^2(A)$ denotes the set of functions in $L^2(\T)$ that vanish outside
%of $A$ and $\widehat f \restrict S$ denotes the restriction of $\widehat f$ to
%$S$.)
%
%\procl c.setdual
%Let $A \subset \T$ be Lebesgue measurable with measure $|A|$.
%Then there is a sequence of density $|A|$ in $\Z$ that is a restriction set
%for $A$.
%Indeed,
%let $\P^A$ be the determinantal probability measure on $2^\Z$
%corresponding to the Toeplitz matrix $(j, k) \mapsto \widehat{\I{A}}(k-j)$.
%Then $\P^A$-a.e.\ $S \subset \Z$ has 
%density $|A|$ in $\Z$ and is a restriction set for $A$.
%\endprocl
%
%Say that a set $S \subseteq \Z$ is {\bf generic} if both
%$$
%{1 \over N} \sum_{0 \le n < N}\delta_{(S-n)} 
%\qquad\hbox{ and } \qquad
%{1 \over N} \sum_{0 \le n < N}\delta_{(-S-n)} 
%$$
%converge weak${}^*$ and to the same limit measure on $2^\Z$. Here,
%$\delta_W$ denotes the unit point mass at $W \in 2^\Z$.
%Every ergodic measure $\mu$ on $2^\Z$ is concentrated on generic sets by the
%ergodic theorem; in this case, the above limits are simply $\mu$.
%Thus, we can strengthen \ref c.setdual/ to say that $S$ is a generic set.
%
%An application of this is the following.
If $A \subseteq \T$ is measurable, recall that $S \subseteq \Z$ is {\bf
complete for} $A$ if the set $ \{ f \I A \st f \in L^2(\T),\, \widehat f
\restrict (\Z \setminus S) \equiv 0 \}$ is dense in $L^2(A)$. 
(Again, $L^2(A)$ denotes the set of functions in $L^2(\T)$ that vanish outside
of $A$ and $\widehat f \restrict S$ denotes the restriction of $\widehat f$ to
$S$.)
%(This is the
%concept dual to restriction set, with the roles of $A$ and $S$ reversed. It is
%also the same as a uniqueness set for $A$, meaning that if $f \in L^2(A)$
%satisfies $\widehat f \restrict S \equiv 0$, then $f \equiv 0$.)
The case where $A$ is an interval is quite classical; see,
e.g., \ref b.Redheffer/ for a review.
A crucial role in that case is played by the following notion of density of
$S$.

\procl d.BM
For an interval $[a, b] \subset \Z \setminus \{ 0 \}$, define its {\bf aspect} 
$$
\asp([a, b])
:=
\max \{ |a|, |b| \}/ \min \{ |a|, |b| \}
\,.
$$
For $S \subseteq \Z$, the {\bf Beurling-Malliavin density} of $S$, denoted
$\bm(S)$, is the supremum of those $D \ge 0$ for which there exist disjoint
nonempty intervals $I_n \subset \Z \setminus \{ 0 \}$ 
with $|S \cap I_n| \ge D |I_n|$ for all $n$ and $\sum_{n \ge 1}
[\asp(I_n)-1]^2 = \infty$.
\endprocl

A simpler form of the Beurling-Malliavin density was provided by \ref
b.Red:Two/, who showed that 
$$
\bm(S) = \inf \left\{ c \st \exists \hbox{ an injection }\beta : S \to \Z
\hbox{ with } \sum_{k \in S} \Big|{1 \over k} - {c \over \beta(k)}\Big| <\infty
\right\} 
\,.
\label e.BMRed
$$

\procl c.seqdual
Let $A \subset \T$ be Lebesgue measurable with measure $|A|$.
Then there is a set of Beurling-Malliavin density $|A|$ in $\Z$ that is 
complete for $A$.
Indeed,
let $\P^A$ be the determinantal probability measure on $2^\Z$
corresponding to the Toeplitz matrix $(j, k) \mapsto \widehat{\I{A}}(k-j)$.
Then $\P^A$-a.e.\ $S \subset \Z$ is complete for $A$ and has
$\bm(S) = |A|$.
\endprocl

When $A$ is an interval, the celebrated theorem of \ref b.BM/ says that if
$S$ is complete for $A$, then $\bm(S) \ge |A|$. (This holds for $S$ that
are not necessarily sets of integers, but we are concerned only with $S
\subseteq \Z$.) 
A.~Ulanovskii has pointed out to the author that this inequality can fail
dramatically for {\it certain\/} subsets $A \subset \T$. We wonder how
small $\bm(S)$ can be for general $A$ and $S \subseteq \Z$ that is complete
for $A$.

\proof
We shall apply \ref t.basis/ with $E := \Z$.
We use the Fourier isomorphism $f \mapsto \widehat f$ between
$L^2(\T)$ and $\ell^2(\Z)$. % and the natural action of $\Z$ on $\ell^2(\Z)$.
Let $H$ be the image of $L^2(A)$ under this isomorphism.
%then $H$ is $\Z$-invariant.
Calculation shows that $\P^A = \P^H$ and that $\P^A[e \in S] = |A|$ for all $e
\in \Z$ when $S$ has law $\P^A$.
%Since $H$ is $\Z$-invariant, $\P^A$ is $\Z$-invariant.
%In addition, $\P^A$ is mixing (whence ergodic) since for every $e \in \Z$, we
%have that $\sum_{f \in \Z} |\ip{P_H e, f}|^2 = \|P_H e\|^2 <\infty$, which
%shows that $\ip{P_H e, f} \to 0$ as $f \to\infty$.
%%(One could alternatively use the Kochen-Stone version of the Borel-Cantelli
%%lemma, which applies because of negative correlations.)
%Ergodicity tells us that
%the density of $S$ is $\P^A$-a.s.\ equal to $\P^A[0 \in S] =
%\widehat{\I{A}}(0) = |A|$.
The equation $\overline{[P_H S]} = H$ is precisely
the statement that $S$ is complete for $A$.
Thus, \ref t.basis/ tells us that $\P^A$-a.e.\ $S$ is complete for $A$.
It remains to show that $\bm(S) = |A|$ $\P^A$-a.s.
If the events $\{ e \in S \}$ were independent for $e \in \Z$, this would be a
special case of a theorem of \ref b.SeipUl/. It is easy to check that
the negative association of $\P^A$ (\ref t.FM/) allows the proof of
\ref b.SeipUl/ to carry through to our situation, just as it does for \ref
c.subG/. 
%(Indeed, the proof of \ref b.SeipUl/ could be based on \ref e.subG/.)
In fact, here is a much shorter proof.
First, the ergodic theorem guarantees that the ordinary density of $S$ is
$|A|$ for $\P^A$-a.e.\ $S$. Thus, $\bm(S) \ge |A|$ $\P^A$-a.s.
For the converse inequality, it suffices by symmetry to consider $S \cap
\Z^+$, which we write as the increasing sequence $\Seq{ s_n \st n \ge 1}$.
By \ref c.subG/ and the Borel-Cantelli lemma, we have 
$$
\bigg| \big|S \cap [1, k]\big| - \big|A\big| k \bigg| \le \sqrt{k \log k}
$$
for all but finitely many $k$ a.s.
If we substitute $k := s_n$, then we obtain that 
$$
\bigg| {1 \over s_n} - {|A| \over n} \bigg|
\le
{\sqrt{s_n \log s_n} \over n s_n}
=
{\sqrt{\log s_n} \over n \sqrt{s_n}}
$$
for all but finitely many $n$ a.s.
Since $s_n \sim n/|A|$ a.s., it follows that 
$$
\sum_n \Big| {1 \over s_n} - {|A| \over n} \Big| <\infty
\quad \hbox{a.s.}
$$
Thus, $S \cap \Z^+$ a.s.\ satisfies \ref e.BMRed/ with $\beta(s_n) := n$.
\Qed
%Since the ordinary density of $S$ exists and equals $|A|$ $\P^A$-a.s.\ by
%the ergodic theorem, and since $\bm(S)$ is always at least the ordinary
%density of $S$, we have left to show that $\bm(S) \le |A|$ $\P^A$-a.s
%In other words, if $D > |A|$, then $\P^A$-a.s.\ there do not exist disjoint
%intervals $I_n$ that satisfy $|S \cap I_n| \ge D |I_n|$ for all $n$ and
%$\sum_{n \ge 1} [\asp(I_n)-1]^2 = \infty$.

The dual form of \ref c.seqdual/ that results from \ref t.morris/ is
that the restriction of the Fourier coefficients of $f$ to $S$ for $f \in
L^2(A)$ is dense in $\ell^2(S)$ for $\P^A$-a.e.\ $S$.
This is an equivalent form of \ref c.not0/.

Two more applications of \ref t.basis/ give results for the wired and free
spanning forests. However, we do not know their significance.
Possibly the one for the $\fsf$, namely, 
$$
\overline{[P_\CYCLE^\perp \fo]} = \CYCLE^\perp  \quad \fsf\hbox{-a.s.}
\,,
$$
could be used to glean more information about the $\fsf$, since it is a
statement about how large $\fo$ must be.

%The ergodicity property that we have used here has a
%much stronger and wider extension, namely, tail triviality of all measures
%$\P^H$. We define and prove this next.
We next define and prove tail triviality of all measures $\P^H$.
For a set $K
\subseteq E$, recall that $\F(K)$ denotes the $\sigma$-field of events
that are measurable with respect to $K$. 
Define the {\bf tail} $\sigma$-field to be the
intersection of $\F(E\setminu K)$ over all finite $K$.
We say that a measure $\P$ on $2^E$ has {\bf trivial tail} if 
every event in the tail $\sigma$-field has measure either 0 or 1.
Recall that tail triviality is equivalent to
$$
\forall \A_1 \in \F(E) \;\ \forall \epsilon > 0 \;\ \exists K \hbox{
finite } \;
\forall \ev A_2 \in \F(E\setminu K)\qquad
\bigl|\P(\A_1 \cap \ev A_2) - \P(\A_1) \P(\ev A_2)\bigr| < \epsilon\,.
\label e.georgii
$$
(See, e.g., \ref b.Georgii:book/, p.~120.)

\procl t.tail
The measure $\P^H$ has trivial tail.
\endprocl

Our proof is modelled on the quantitative proof of tail triviality of
$\fsf$ and $\wsf$ in \BLPSusf.
We explain what changes 
are needed to make that proof work here (and, along the way, give
slight corrections).
The quantitative form of tail triviality that we prove is this:

\procl t.correl
Let $F$ and $K$
be disjoint nonempty subsets of $E$ with $K$ finite.
Let $C$ be a subset of $K$.  
Then 
$$
\Var^H\big(\P^H[\ba \cap K = C\mid \F(F)]\big) \le |K| \sum_{e \in K}
\|P_{[F]} P_H(e) \|^2
\,.
\label e.correl
$$
If $\A_1 \in \F(K)$ and $\A_2 \in \F(F)$, then
$$
|\P^H(\A_1 \cap \A_2) - \P^H(\A_1) \P^H(\A_2)| 
\le
\left(2^{|K|} |K| \sum_{e \in K} \|P_{[F]} P_H(e) \|^2\right)^{1/2}
\,.
\label e.correl2
$$
\endprocl

Before proving \ref t.correl/, we explain why it implies
\ref  t.tail/. In fact, we show the more quantitative \ref e.georgii/.
Let $\A$ be any event and $\epsilon > 0$. Find a finite set $K_1$ and
$\A_1 \in \F(K_1)$ such that $\P^H(\A_1 \triangle \A) < \epsilon /3$.
Now find a finite set $K_2$ so that
$$
\left(2^{|K_1|} 
|K_1| \sum_{e \in K_1} \|P_{[E \setminus K_2]} P_H(e) \|^2\right)^{1/2}
< \epsilon /3\,.
$$
Then for all $\A_2 \in \F(E \setminu K_2)$, we have $|\P^H(\A \cap \A_2) -
\P^H(\A)\P^H(\A_2)| < \epsilon$.

To prove \ref t.correl/, we need to establish some lemmas.
Note that both sides of \ref e.correl/, as well as of \ref e.correl2/, are
continuous for an increasing sequence of sets $F$, whence it suffices to
prove both inequalities only for $F$ finite.
Thus, we may actually assume that $E$ is finite.
Assume now that $E$ is finite.
Let $Q_F$ be the orthogonal projection onto the (random) subspace $H^F_{F\cap
\ba}$ defined in \ref e.HFSdef/.  

\procl l.conds
Let 
$F\subset E$.  Then
$$
\E^H Q_F
=
\sum_{S\subseteq F} \P^H[ \ba \cap F=S] P_\HFS
=
P_{[F]}^\perp P_H P_{[F]}^\perp 
\,.
\label e.PHP
$$
\endprocl

The proof is the same as that of Lemma 8.5 of \BLPSusf, where we now
establish
that for any $e, e' \in E$, we have
$$
\E^H \ip{Q_F e, e'} 
=
\ip{P_{[F]}^\perp P_H P_{[F]}^\perp e, e'} 
\,.
%\label e.PHPeh
$$
This uses
\ref p.kirchproj/ in place of the direct arguments in \BLPSusf.

The next lemma has a precisely parallel proof to that of Lemma 8.6 of
\BLPSusf.

\procl l.oneentry
Let $F\subset E$
and $u\in\H$.
Then
$$
\Var^H(Q_F u)
:=
\E^H\Big[\| Q_F u - \E^H Q_F u \|^2\Big]
= \|P_{[F]} P_H P_{[F]}^\perp u\|^2
\,.
$$
\endprocl

% It is interesting to note that similar
% calculations show that if $Q$ is a 
% random projection in a Hilbert space, then
% $$
% \Var(Q\xi) = (\E Q\xi,\E Q^\perp \xi)
% \,.
% $$

\proofof t.correl
Define $\QF$ to be the orthogonal projection on the subspace $H_{F
\cap B, F \setminus B}$.
By \ref t.genprs/ and \ref e.HABcond/, we have
$$
\P^H[ \ba \cap K = C\mid  \F(F)](B)
=
\det \CC_{C, B}^K
\,,
$$
where
$$
\CC_{C, B}^K
:=
\Bigl[\bigip{\QF^{C,e}  e,  {e'}}\Bigr]_{e, e' \in K}
$$
with notation as follows:
For a set $C$ and operator $P$, write
$$
P^{C, e} := \cases{P        &if $e \in C$,\cr
                   \id-P    &if $e \notin C$.\cr}
\label e.Be
$$
Since $K \cap F = \emptyset$, we have that $\bigip{\QF^{C,e}  e,  {e'}} =
\bigip{Q_F^{C,e}  e,  {e'}}$ for $e, e' \in K$ on the event that $\ba = B$.
Thus
\begineqalno
\E^H \CC_{C, \ba}^K 
&=
[\ip{\E^H Q_F^{C,e}  e,  {e'}}]_{e, e' \in K}
\cr&=
\Bigl[\bigip{(P_{[F]}^\perp P_H P_{[F]}^\perp)^{C,e}   e,
           {e'}}\Bigr]_{e, e' \in K}
\qquad\qquad\hbox{by \ref e.PHP/}
%\cr&=
%\Bigl[\bigip{P_{[F]}^\perp P_H^{C,e} P_{[F]}^\perp   e,
           %{e'}}\Bigr]_{e, e' \in K}
\cr&=
\Bigl[\bigip{P_H^{C,e}  e,  {e'}}\Bigr]_{e, e' \in K}
\qquad\qquad\hbox{since $K \cap F = \emptyset$}
.
\cr
\endeqalno
Therefore, 
$$
\E^H\det\CC_{C, \ba}^K 
=
\E^H\P^H[ \ba \cap K = C\mid  \F(F)]
=
\P^H[ \ba \cap K = C]
=
\det\E^H\CC_{C, \ba}^K
$$
by \ref t.genprs/.
Furthermore, for any orthogonal projection $P$, we have
$$
\sum_{e' \in K} |\ip{P  e,  {e'}}|^2
\le \|P  e \|^2
\le 1
$$
because $\Seq{ {e'}\st e'\in E}$ is an orthonormal basis for $\H$.
Thus, we may apply Lemma 8.7 of \BLPSusf\ to obtain
\begineqalno
\Var^H\big(\P^H[ \ba &\cap K = C\mid  \F(F)]\big)
=
\Var^H\left(\det \CC_{C, \ba}^K\right)
\cr&\le
|K| \sum_{e, e' \in K} \Var^H \bigip{Q_F^{C,e}  e,  {e'}}
\le
|K| \sum_{e\in K, e' \in E} \Var^H \bigip{Q_F^{C,e}  e, {e'}}
\cr&=
|K| \sum_{e\in K} \Var^H \bigl(Q_F^{C,e} e\bigr) 
=
|K| \sum_{e\in K} \Var^H \bigl(Q_F e\bigr) 
\cr&=
|K| \sum_{e \in K} \|P_{[F]} P_H(e)\|^2
\,,
\cr
\endeqalno
using \ref l.oneentry/.
This proves \ref e.correl/.

To deduce \ref e.correl2/ from \ref e.correl/, write $a := 2^{|K|}
|K| \sum_{e \in K} \|P_{[F]} P_H(e) \|^2$. Then for all $\A_1 \in \F(K)$,
we have
$$
\Var^H\big(\P( \A_1 \mid  \F(F))\big) \le a
$$
since $\A_1$ is the union of at most $2^{|K|}$ disjoint
cylinder events of the form $ \{ \ba \cap K = C \}$.
Therefore for all $\A_2 \in \F(F)$,
$$
\big|\P^H(\A_1 \mid \A_2) - \P^H(\A_1)\big|^2\, \P^H(\A_2) \le a
\,,
$$
so that
$$
\big|\P^H(\A_1 \cap \A_2) - \P^H(\A_1) \P^H(\A_2)\big|^2 \le a \P^H(\A_2) \le a
\,.
$$
This is the same as \ref e.correl2/.
\Qed

\bsection{Positive Contractions}{s.contract}

We have seen that the matrix of any orthogonal projection gives a
determinantal probability measure.
We now do the same for positive contractions and give their properties.

We call $Q$ a {\bf positive contraction} if $Q$ is a self-adjoint
operator on $\H$ such that for all $u \in \H$, we have $0 \le (Q u, u) \le
(u, u)$.
To show existence of a corresponding determinantal probability measure, which
we shall denote $\P^Q$, let $P_H$ be any orthogonal projection that is
a {\bf dilation} of $Q$, i.e., $H$ is
a closed subspace of $\ell^2(E')$ for some $E' \supseteq E$ and for all $u
\in \H$, we have $Qu = P_{\H} P_H u$, where we regard $\ell^2(E')$ as the
orthogonal sum $\H \oplus \ell^2(E' \setminus E)$. (In this case, $Q$
is also called the {\bf compression} of $P_H$ to $\H$.)
The existence of a dilation is standard and is easily
constructed: Let $E'$ be the union of $E$ with a disjoint copy $\hat E$ of
$E$. Let $T$ be the positive square root of $Q$  and let $\hat T$
be the positive square root of $I-Q$. 
The operator whose block matrix is 
$$
\left(\matrix{Q&T \hat T\cr T \hat T&I-Q\cr}\right)
$$
is easily checked to be self-adjoint and idempotent, hence it is an orthogonal
projection onto a closed subspace $H$.
Having chosen a dilation, we simply define $\P^Q$ as the law of $\ba \cap
E$ when $\ba$ has the law $\P^H$.
Then \ref e.DPM/ is a special case of \ref e.included/.

Of course, when $Q$ is the orthogonal projection onto a subspace $H$, then
$\P^Q = \P^H$.

The basic properties of $\P^Q$ follow from those for orthogonal projections.
In the following, we write $Q_1 \le Q_2$ if $\ip{Q_1 u, u} \le \ip{Q_2 u, u}$
for all $u \in \H$.

\procl t.Q
Let $Q$ be a positive contraction.
For any finite $A, B \subseteq E$, we have 
$$
\P^Q\left[ A \subseteq \qba, B \cap \qba = \emptyset\right]
=
\leftip{\bigwedge_{e \in A} Q e \wedge \bigwedge_{e \in B} (I-Q) e,
\theta_A \wedge \theta_B }
\,.
\label e.Qgenprs
$$
The measure $\P^Q$ has conditional negative associations (with external
fields) and a trivial tail $\sigma$-field.
If $A \subseteq E$ is finite and $\mu := \E^Q[|\qba \cap A|]$, then for any
$a > 0$, we have 
$$
\P^Q\bigg[\Big| |\qba \cap A| - \mu \Big|  \ge a\bigg] \le 2 e^{-2 a^2/|A|}
\,.
\label e.QsubG
$$
If $Q_1$ and $Q_2$ are commuting positive contractions and $Q_1 \le Q_2$,
then $\P^{Q_1} \preccurlyeq \P^{Q_2}$.
\endprocl

\procl r.trivtail
Independently, \ref b.ShiTak:II/ showed that $\P^Q$ has a trivial tail
$\sigma$-field when the spectrum of $Q$ lies in $(0, 1)$. 
\endprocl

\proof
The first four properties are immediate consequences of \ref e.genprs/,
\ref t.FM/, \ref t.tail/, and \ref c.subG/.
(Of course, many other properties follow from the negative association; we
mention \ref e.QsubG/ merely as an example.)
The last statement will follow from \ref t.dominate-infinite/ once we show
that the hypothesized commutativity implies that we may take dilations
$P_{H_i}$ of $Q_i$ with $H_1 \subseteq H_2$.
To do this, we use the 
following form of the spectral theorem:
There is a measure space $(X, \mu)$, two Borel
functions $f_i : X \to [0, 1]$, and a unitary map
$U : \H \to L^2(\mu)$ such that 
$U Q_i U^{-1} : g \mapsto f_i g$ for $i=1, 2$
and any $g \in L^2(\mu)$ (apply Theorem IX.4.6, p.~272, of \ref
b.Conway:CFA/ to the normal operator $Q_1 + i Q_2$).
Since $Q_1 \le Q_2$, we have $f_1 \le f_2$.
Use $U$ to identify $\H$ with $L^2(\mu)$ and to identify $Q_i$ with $M_i := U
Q_i U^{-1}$.
Let $\lambda$ denote Lebesgue measure on $[0, 1]$ and define $A_i := \{
(x, y) \st x \in X,\, 0 \le y \le f_i(x) \}$.
Let $H_i$ be the subspace of functions in $L^2(\mu \otimes
\lambda)$ that vanish outside $A_i$.
Since $A_1 \subseteq A_2$, we have $H_1 \subseteq H_2$.
Embed $L^2(\mu)$ in $L^2(\mu \otimes \lambda)$ by
$g \mapsto g \otimes \constant 1$ and identify
$L^2(\mu)$ with its image.
%Extend $U$ to $\widetilde U : \H \oplus \bigoplus_j L^2(\lambda) \to
%\bigoplus_j L^2(\mu_j \otimes \lambda)$ by $\widetilde U : (v, \Seq{h_j})
%\maptso \Seq{(U v)_j \otimes h_j}$.
Then $M_i$ is the compression of $P_{H_i}$ to $L^2(\mu)$, as
desired.
\comment{Could also use resolution of identity.}
\Qed

A formula for the probability measure $\P^Q(\;\cbuldot \mid A \subseteq \qba,
B \cap \qba = \emptyset)$ follows from applying \ref e.HABcond/ to a dilation
of $Q$.
However, this is not very explicit. Often conditioning on just $A \subseteq
\qba$ is important, so we give the following direct formula for that case.
Note that we allow $A$ to be infinite; if $A = \bigcup_n A_n$ with $A_n$
finite, then $\P^Q(\;\cbuldot \mid A_n \subseteq \qba)$ is a stochastically
decreasing sequence of probability measures by \ref t.FM/ and so defines
$\P^Q(\;\cbuldot \mid A \subseteq \qba)$.
We shall write 
$$
\ip{u, v}_Q := \ip{Qu, v}
$$
for the inner product on $\H$ induced by $Q$.
Let $[E]_Q$ be the completion of $\H$ in this inner product
and $P^\perp_{[A]_Q}$ be the orthogonal projection in $[E]_Q$ onto the
subspace orthogonal to $A$.

\procl p.condonA
Let $Q$ be a positive contraction on $\H$ and $A \subset E$.
When $\qba$ has law $\P^Q$ conditioned on $A \subseteq \qba$, then the law of
$\qba \cap (E \setminus A)$ is the determinantal probability measure
corresponding to the positive contraction on $\ell^2(E \setminus A)$ whose
$(e, f)$-matrix entry is 
$$
\bigip{P^\perp_{[A]_Q} e, P^\perp_{[A]_Q} f}_Q
\,.
$$
\endprocl

An equivalent expression was found independently by \ref b.ShiTak:I/,
Corollary 6.5.

\proof
Because of \ref e.HABcond/, we know that the law of $\qba \cap (E \setminus
A)$ is the determinantal probability measure corresponding to the
compression of some orthogonal projection, i.e., to some positive
contraction.
What remains is to show that for any finite $B \subset E \setminus A$, we have
$$
\P^Q[B \subseteq \qba \mid A \subseteq \qba]
=
\det\big[\bigip{P^\perp_{[A]_Q} e, P^\perp_{[A]_Q} f}_Q\big]_{e, f \in B}
\,.
$$
By continuity, it suffices to do this when $A$ is finite.
Now, 
$$
\P^Q[B \subseteq \qba \mid A \subseteq \qba]
=
{\det[\ip{Q e, f}]_{e, f \in A \cup B}
\over
\det[\ip{Q e, f}]_{e, f \in A}}
=
{\det[\ip{e, f}_Q]_{e, f \in A \cup B}
\over
\det[\ip{e, f}_Q]_{e, f \in A}}
=
{\|\theta_B \wedge \theta_A\|_Q^2 \over \|\theta_A\|_Q^2}
\,.
$$
We use the following fact about exterior algebras.
For any vectors $u_1, u_2, \ldots, u_m, v_1, v_2, \ldots, v_n$ with $H$
defined to be the span of $v_1, \ldots, v_n$, we have 
$$
\bigwedge_{i=1}^m u_i \wedge \bigwedge_{j=1}^n v_j
=
\bigwedge_{i=1}^m P^\perp_H u_i \wedge \bigwedge_{j=1}^n v_j
\,.
$$
This is because $P_H u_i \wedge \bigwedge_{j=1}^n v_j = 0$.
Thus,
$$
\P^Q[B \subseteq \qba \mid A \subseteq \qba]
=
{\|\bigwedge_{e \in B} P^\perp_{[A]_Q} e \wedge \theta_A\|_Q^2
\over \|\theta_A\|_Q^2}
=
\|\bigwedge_{e \in B} P^\perp_{[A]_Q} e\|_Q^2
\,,
$$
as desired.
\Qed

%As we shall see in \ref s.applic:dyn/, there are many quite interesting
%As illustrated by many examples in \ref b.LS:dyn/, e.g., there are numerous
%interesting measures $\P^Q$ for $Q$ not a projection. 

\procl r.incexc
If \ref e.DPM/ is given, then \ref e.Qgenprs/ can be deduced from
\ref e.DPM/ without using our general theory and, in fact, without assuming
that the matrix $Q$ is self-adjoint. Indeed,
suppose that $X$ is any diagonal matrix. Denote its $(e, e)$-entry by $X_e$.
Comparing coefficients of $X_e$ shows that \ref e.DPM/ implies 
$$
\EBig{\prod_{e \in A} \big(\II{e \in \qba} + X_e\big)}
=
\det \big( (Q + X) \restrict A \big)
\,.
$$
If $A$ is partitioned as $A_1 \cup A_2$ and we choose $X$ so that $X_e =
-\I{A_2}(e)$, then we obtain 
$$
\P[A_1 \subseteq \qba, A_2 \cap \qba = \emptyset]
=
\det \big( Q^{A_2} \restrict (A_1 \cup A_2) \big)
\,,
\label e.incexc
$$
where $Q^A$ denotes the matrix whose rows are the same as those of $Q$ except
for those rows indexed by $e \in A$, which instead equals that row of $Q$
subtracted from the corresponding row of the identity matrix. In other words,
if $Q_{e,f}$ denotes the $(e, f)$-entry of $Q$, then the $(e, f)$-entry of
$Q^A$ is 
$$
\I{A}(e) + (-1)^{\I{A}(e)} Q_{e,f}
\,.
$$
This is another form of \ref e.Qgenprs/.
(Equation \ref e.incexc/ amounts to an explicit form of the
inclusion-exclusion principle for a determinantal probability measure.)
\endprocl

\procl r.contract
If $Q$ is a self-adjoint matrix such that \ref e.DPM/ defines a probability
measure, then necessarily $Q$ is a positive contraction.
The fact that $Q \ge 0$ is a consequence of having nonnegative minors,
while $I - Q \ge 0$ follows from observing that $I - Q$ also defines a
determinantal probability measure, the dual to the one defined by $Q$.
\endprocl

%As one would expect, the matroid whose independent sets are the
%intersections of the independent sets of $\M = (E', \B)$ with $E$
%is called the {\bf restriction} of $\M$ to $E$.

\bsection{Open Questions: General Theory}{s.applic:thy}

Our last sections present open questions organized by topic.
The first two sections concern general determinantal probability measures,
while the others examine specific types of measures.

In order to extend \ref e.HABcond/ to the case where $A$ and $B$ may be
infinite and thereby obtain a version of conditional probabilities,
define
$$
H^*_{A, B} := \overline{ ( H \cap A^\perp ) + [A \cup B] } \cap B^\perp
$$
and
$$
H^{**}_{A, B} = \big( \overline{ H + [B] } \cap (A \cup B)^\perp \big) + [A]
\,.
$$
One can show that $H^*_{A, B} \subseteq H^{**}_{A, B}$, but that they are not
necessarily equal.

The following conjecture would greatly simplify the proof of \ref t.morris/
above.

\procl g.version
Let $H$ be a closed subspace of $\H$ and $K \subset E$.
A version of the conditional probability measure $\P^H$ given $\F(K)$ is
$B \mapsto \P^{H^*_{K \cap B, K \setminus B}}$ and another is given by
$B \mapsto \P^{H^{**}_{K \cap B, K \setminus B}}$.
\endprocl

\comment{
\proofof t.morris
%Since $H_{\ba} \subseteq [\ba]$, the event
%$\{e \in H_{\ba} \}$ is contained in the event $\{ e \in \ba \}$.
%On the other hand, if $e \in B$ yet $e \notin H_B$, then by \ref l.HB/, we
Suppose that $e \in B$ yet $e \notin H_B$. Then by \ref l.HB/, we
have $e \notin \overline{H+B^\perp} = \left(H^\perp \cap [B]\right)^\perp$.
Therefore $P_{H^\perp \cap [B]} e \ne 0$, and so there is some $u \in H^\perp
\cap [B]$ with $\ip{u, e} \ne 0$. 
We may write $u = a e + v$ for some $a \ne 0$ and some $v \in e^\perp$, whence
$e \in (H^\perp \cap [B]) + e^\perp$.
Since $e \in (B \setminus \{ e \})^\perp$, it follows that 
$$
e \in (H^\perp)_{E \setminus B, B \setminus \{ e \}} 
\subseteq (H^\perp)^*_{E \setminus B, B \setminus \{ e \}}
= \left( H^*_{B \setminus \{ e \}, E \setminus B} \right)^\perp
\,.
$$
Consequently, $\P^{H^*_{B \setminus \{ e \}, E \setminus B}}[e \in \ba] = 0$.

By \ref t.version/, this means that the $\P^H$-probability that $e \in \ba$
conditioned on $\ba \restrict (E \setminus \{ e \})$ is 0
a.s.\ on the event $\{ e \notin H_{\ba}\}$.
Therefore, $\P^H[e \in \ba, e \notin H_{\ba}] = 0$.
As this holds for all $e$ belonging to the countable set $E$, we have $\ba
\subset H_{\ba}$ $\P^H$-a.s., which is equivalent to the desired
result.
\Qed
}

We say that $\ev A_1,\ev A_2\subset 2^E$ {\bf occur disjointly for}
$F\subseteq E$ if there are disjoint sets $F_1,F_2\subset E$ 
such that 
$$
\{ K \subseteq E \st K\cap F_i=F\cap F_i \} \subseteq \ev A_i
$$
for $i = 1, 2$.
A probability measure $\P$ on $2^E$ is said to have the {\bf BK property}
if 
$$
\P[\ev A_1 \hbox{ and } \ev A_2 \hbox{ occur disjointly for }
\qba ]
\le
\P[\qba \in \ev A_1]\P[\qba \in \ev A_2]
$$
for every pair $\ev A_1, \ev A_2\subset 2^E$ of increasing events.
Does every determinantal probability measure $\P^Q$ have the BK property? 
The BK inequality of \ref b.BK/ says
that this holds when $Q$ is a diagonal matrix, i.e., when $\P$
is product measure.
The answer is unknown even in the special case of uniform spanning trees,
where it is conjectured to hold in
\BLPSusf. 

Is entropy concave in $Q$ for fixed $E$?
That is, for finite $E$ and a positive contraction $Q$, define the {\bf
entropy} of $\P^Q$ to be
$$
\ent(Q) := - \sum_{A \in 2^E} \P^Q[A] \log \P^Q[A] 
\,.
$$
Numerical calculation supports the following conjecture. 

\procl g.concave
For any positive contractions $Q_1$ and $Q_2$, we have 
$$
\ent\big((Q_1+Q_2)/2\big) \ge \big(\ent(Q_1) + \ent(Q_2)\big)/2
\,.
\label e.concave
$$
%Moreover, given positive contractions $Q_1 \ne Q_2$, define the function
%$$
%h(t) := \ent\big( (1-t) Q_1 + t Q_2 \big)
%$$
%for $0 \le t \le 1$.
%For any even integer $n \ge 2$ and $0 < t < 1$, the $n$th derivative
%$f^{(n)}(t) < 0$.
\endprocl

In \BLPSusf, it is asked whether the free and wired
spanning forests are mutually singular when they are not equal.
One might hope that the following more general
statement holds: if $H_1 \subsetneq H_2$, then the corresponding
probability measures $\P^{H_1}$ and $\P^{H_2}$ are mutually singular.
However, this general statement is false, as shown by \ref b.HL:change/.
Nevertheless, since it is true trivially for finite $E$, it seems likely that
there are interesting sufficient conditions for $\P^{H_1}$ and $\P^{H_2}$ to
be mutually singular.

If $\P^{H_1} \preccurlyeq \P^{H_2}$, must there exist a subspace $H_3 \subseteq
H_2$ such that $\P^{H_1} = \P^{H_3}$? (This was answered in the negative by
Lewis Bowen after a preprint was circulated.)

Given the value of \ref t.basis/ and of its dual form \ref t.morris/,
it seems desirable to extend other properties of matroids to the infinite
setting.

\bsection{Open Questions: Coupling}{s.applic:coupling}

A {\bf coupling} of two probability measures $\P^1$, $\P^2$ on
$2^E$ is a probability measure $\mu$
on $2^E \times 2^E$ whose coordinate projections are $\P^1$, $\P^2$, meaning
that for all events $\A \subseteq 2^E$, we have
$$
\mu\big\{ (A_1, A_2) \st A_1 \in \A \big \}
=
\P^1(\A)
$$
and
$$
\mu\big \{(A_1, A_2) \st A_2 \in \A \big \}
=
\P^2(\A)
\,.
$$
A coupling $\mu$ is called
{\bf monotone} if 
$$
\mu\big\{(A_1, A_2) \st A_1 \subseteq A_2\big\} = 1
\,.
$$
By \inviscite{Strassen}Strassen's \refbyear{Strassen} theorem, stochastic
domination $\P^1 \preccurlyeq \P^2$ is equivalent to the existence of a
monotone coupling of $\P^1$ and $\P^2$.
A very interesting open question that arises from \ref
t.dominate/ is to find a natural or explicit monotone coupling of $\P^{H'}$ and
$\P^{H}$ when $H' \subset H$.  

A coupling $\mu$ is {\bf disjoint} 
if $\mu\big\{(A_1, A_2) \st A_1 \cap A_2 = \emptyset\big\} = 1$.
A coupling $\mu$ has 
{\bf union marginal} $\P$ if
for all events $\A \subseteq 2^E$, we have
$$
\P(\A) = \mu\big\{(A_1, A_2) \st A_1 \cup A_2 \in \A \big \}
\,.
$$

\procl q.unioncoupling
Given $H = H_1 \oplus H_2$, is there a (natural or otherwise) disjoint
coupling of $\P^{H_1}$ and $\P^{H_2}$ with union marginal $\P^H$?  
\endprocl

This is easily seen to be the case when
$H=\H$: The probability measure $\mu$ on $2^E \times
2^E$ defined by 
$$
\mu\{(A, E \setminu A) \st A \in \A\} := \P^{H_1}(\A)
$$
and 
$$
\mu\big \{ (A, B) \st B \ne E \setminus A \} := 0
$$
does this, as we can see by \ref c.dualrep/.  
A positive answer in general to \ref q.unioncoupling/
would give the following more general
result (by the method of proof of \ref t.Q/): If $Q_1$ and $Q_2$ are
commuting positive contractions on $\H$ such that $Q_1 + Q_2 \le I$, then
there is a disjoint coupling of $\P^{Q_i}$ with union marginal $\P^{Q_1+Q_2}$.
We note that the requirement of being disjoint is superfluous, although useful
to keep in mind:

\procl p.disjoint
If $Q_1$ and $Q_2$ are positive contractions on $\H$ such that $Q_1
+ Q_2 \le I$, then any coupling of $\P^{Q_1}$, $\P^{Q_2}$ with union marginal
$\P^{Q_1+Q_2}$ is necessarily a disjoint coupling.
\endprocl

\proof 
Let the coupling be $\mu$, which picks a random pair $(\qba_1, \qba_2) \in 2^E
\times 2^E$.
Then for all $e \in E$, we have
\begineqalno
\bigip{(Q_1+Q_2)e, e}
&=
\P^{Q_1+Q_2}[e \in \qba]
=
\mu[e \in \qba_1 \cup \qba_2]
\cr&=
\mu[e \in \qba_1] + \mu[e \in \qba_2] - \mu[e \in \qba_1 \cap \qba_2]
\cr&=
\P^{Q_1}[e \in \qba_1] + \P^{Q_2}[e \in \qba_2]- \mu[e \in \qba_1 \cap \qba_2]
\cr&=
\ip{Q_1 e, e} + \ip{Q_2 e, e}- \mu[e \in \qba_1 \cap \qba_2]
\cr&=
\bigip{(Q_1+Q_2)e, e}- \mu[e \in \qba_1 \cap \qba_2]
\,.
\endeqalno
Therefore $\mu[e \in \qba_1 \cap \qba_2] = 0$.
Since this holds for each $e$, we get the result.
\Qed

\comment{
It is not true that a coupling can be arrived at via a simple argument with
a projection measure on $\R^{\edge \cup \edge'}$, since numbers of elements of
$\edge$ and of $\edge'$ would be constant, yet negatively correlated. This
means they are independent, but then they cannot be disjoint.
}

If a natural monotone coupling is found, it ought to provide a coupling of
the free and wired spanning forests that is invariant under all
automorphisms of the underlying graph, $G$. This should help in
understanding the free spanning forest.
(A specific instance is given below.)

We shall give some partial results on the general question.

\procl p.diff
If $H \subset H'$ are two closed subspaces of $\H$, then there is a
monotone coupling of $\P^{H}$ and $\P^{H'}$ concentrated on the set
$\{(B_1, B_2) \st |B_2 \setminu B_1| = k\}$, where $k$ is the codimension
of $H$ in $H'$ (possibly $k = \infty$).
\endprocl

We do not know whether every monotone coupling has this property.

Of course, \ref p.diff/ is trivial when $|E| < \infty$. To prove \ref p.diff/
when $E$ is infinite, we first prove a lemma that shows that in the
finite-dimensional codimension-one case, every monotone coupling gives rise to
a disjoint coupling with the proper union marginal:

\procl l.codim1
Let $H$ be a finite-dimensional subspace of $\H$ and $u$ be a unit vector
in $H^\perp$. Let $H'$ be the span of $H$ and $u$. If $\mu$ is any monotone
coupling of $\P^H$ and $\P^{H'}$, then for every event $\ev A$,
$$
\mu{\{(B, B') \st B' \setminu B \in \ev A\}} = \P^u(\ev A)
\,.
$$
\endprocl

\proof
We have that 
$$
\mu{\{(B, B') \st |B| = \dim H \}} = 1
$$
and
$$
\mu{\{(B, B') \st |B'| = \dim H + 1 \}} = 1
\,.
$$
Therefore, it suffices to show that for every $e \in E$, 
$$
\mu{\{(B, B') \st e \in B' \setminu B \}} = \P^u(\{e\})
\,.
$$
This equation holds because the left-hand side is equal to
\begineqalno
\mu{\{(B, B') \st e \in B'\}} 
-
\mu{\{(B, B') \st e \in B \}}
&=
\P^{H'}[e \in \ba'] - \P^H[e \in \ba]
\cr&=
\|P_{H'}e\|^2 - \|P_H e\|^2
=
\|P_{[u]} e\|^2
\cr&=
\P^u[\{e\}]
\,.
\Qed
\endeqalno

It follows from this and duality considerations that \ref q.unioncoupling/ has
a positive answer whenever $|E| \le 5$.
We have tested by computer thousands of random instances of \ref
q.unioncoupling/ for $|E| = 6, 7, 8, 9$ and all have a positive answer.

\proofof p.diff
The case that $H'$ is finite dimensional is trivial, so suppose that $H'$
is infinite dimensional.

Suppose first that $k =1$. Let $H' = H \smalloplus [u]$ and choose any
increasing sequence of finite-dimensional subspaces $H_i$ ($i \ge 1$) of
$H$ whose union is dense in $H$. Let $\mu_i$ be any monotone coupling of
$\P^{H_i}$ and $\P^{H_i \smalloplus [u]}$ for each $i \ge 1$. Let $\mu$ be
any weak${}^*$ limit point of $\mu_i$. Then $\mu$ is a monotone coupling of
$\P^{H}$ and $\P^{H'}$. Furthermore, since
$$
\mu_i{\{(B, B') \st B' \setminus B = \{e\} \}} = \P^u[\{e\}]
$$
for each $i$, the same holds for $\mu$. This gives the desired conclusion.

Now suppose that $1 < k < \infty$.
Let $H' = H_1 \supset H_2 \supset H_3 \supset \cdots \supset H_k = H$ be a
decreasing sequence of subspaces with $H_{i+1}$ having codimension 1 in
$H_i$ for each $i = 1, \ldots, k-1$. By what we have shown, we may choose
monotone couplings $\mu_i$ of $\P^{H_{i+1}}$ with $\P^{H_i}$ for each $i =
1, \ldots, k-1$ with the property that
$$
\mu_i{\{(B, B') \st |B' \setminus B| = 1 \}} = 1
\,.
$$
Choose $(B_i, B'_i)$ with distribution $\mu_i$ and independently of each
other. Then the distribution $\mu$ of $(B_{k-1}, B'_1)$ given that $B_i =
B'_{i+1}$ for each $i = 1, \ldots, k-1$ is the desired coupling.

Finally, if $k = \infty$, then choose a decreasing sequence of subspaces
$H' = H_1 \supset H_2 \supset H_3 \supset \cdots \supset H_i \supset \cdots$
with
$H_{i+1}$ having codimension 1 in $H_i$ for each $i  \ge 1$ and $\bigcap
H_i = H$. Let $\mu_i$ be
any monotone coupling of $\P^{H'}$ and $\P^{H_i}$ with
$$
\mu_i{\{(B, B') \st |B' \setminus B| = i \}} = 1
\,.
$$
Let $\mu$ be any weak${}^*$ limit point of $\mu_i$. Then $\mu$ is the desired
coupling.
\Qed

An example of the usefulness of coupling is as follows.
Let $G= (\vertex, \edge)$ be a proper planar graph with planar dual $\d G=
(\d\vertex, \d\edge)$. Let $H$ be a subspace of $\ell^2(\edge)$.
Let $\phi$ be the map $e \mapsto \d e$.
%We may identify $\edge$ and $\d\edge$. 
Then $\phi$ induces a map $\ell^2(\edge) \to \ell^2({\d \edge})$ that sends
$H$ to a subspace $\d H$ of $\ell^2({\d \edge})$.
For example, $\d{\STAR(G)} = \CYCLE(\d G)$ and $\d{\CYCLE(G)} =
\STAR(\d G)$. %, and $\d{\grad\HD(G)} = \grad\HD(\d G)$.
The disjoint coupling of $\P^H$ and $\P^{H^\perp}$ of the first paragraph
of this section gives
a coupling of $\P^H$ on $\ell^2(\edge)$ and $\P^{(\d H)^\perp}$ on
$\ell^2({\d\edge})$ for which exactly one of $e$ and $\d e$ appear in $\ba$ and
$\d \ba$ for each $e \in \edge$.
For example, if $H=\STAR(G)$, then $\P^H = \wsf(G)$; since $(\d H)^\perp
= \CYCLE(\d G)^\perp$, we obtain the disjoint coupling of $\wsf(G)$
with $\fsf(\d G)$ used in \BLPSusf.  This is the just about the only
method known to derive much information about $\fsf$ when $\fsf \ne \wsf$.
Note that if $G$ is not planar but is embedded on a surface, then
$\d{\STAR(G)}$ is the span of the facial (contractible) cycles of $\d
G$, which may be smaller than $\CYCLE(\d G)$.
(This is discussed further in \ref s.applic:CW/ below.)

\bsection{Open Questions: Groups}{s.applic:Groups}

We shall consider first finite groups, then infinite groups.

Suppose that $E$ is a finite group.
Then $\H$ is the group algebra of $E$.
Invariant subspaces $H$ give subrepresentations of the regular
representation and give invariant probability measures $\P^H$.
There is a canonical decomposition 
$$
\H = \bigoplus_{j=1}^s  H_j
\,,
$$
where each $H_j$ is an invariant subspace containing all isomorphic
copies of a given irreducible representation.
(See, e.g., \ref b.FultonHarris:book/.)
The matrix of $P_{H_j}$ is given by the character
$\chi_{H_j}$ of the representation, namely, the $(e, f)$-entry is
$\overline{\chi_{H_j}(e f^{-1})}/|E|$
(\ref b.FultonHarris:book/, p.~23, (2.32)).
Can we (disjointly) couple all measures $\P^{H_j}$ so that every
partial union has marginal equal to $\P^H$ for $H$ the
corresponding partial sum?
In other words, is there a probability measure $\mu$ on $\prod_{j=1}^s 2^E$
picking a random $s$-tuple $\Seq{\qba_1, \ldots, \qba_s}$
such that for every $J \subseteq \{1, \ldots, s\}$, the law of $\bigcup_{j \in
J} \qba_j$ is $\P^{H_J}$, where $H_J := \bigoplus_{j \in J}
H_j$?
We call such a coupling {\bf complete}.

Consider the case $E = \Z_n$.
All irreducible representations are 1-dimensional and there are $n$ of them:
for each $k \in \Z_n$, we have the representation 
$$
m \mapsto e^{2\pi i k m/n} \qquad (m \in \Z_n)
\,.
$$
Thus, a complete coupling would be a random permutation of $\Z_n$ with special
properties.
By averaging, we may always assume that any complete coupling is invariant.
Do they always exist? If so, the set of invariant complete couplings is a
polytope. What are
its extreme points or the supports of the extreme points? What is its
barycenter? What asymptotic properties distinguish it from the uniform
permutation?
Testing by computer indicates existence for all $n \le 7$.
Thus, it would appear that complete couplings always exist on $\Z_n$.

One should be aware that
it is not always possible to completely couple 4 measures from orthogonal
subspaces when the subspaces are not invariant, as one can show from the
following example. Let $v_1 := \Seq{1, 1, -3, 1}$, $v_2 := \Seq{1, -1, 5, 2}$,
$v_3 := \Seq{1, 1, -2, -2}$, and $v_4 := \Seq{-3, 2, 1, 4}$.
Let $u_j$ be the corresponding vectors resulting from the Gram-Schmidt
procedure, that is, $u_1 := v_1/\|v_1\|$, $u_2 := P_{[v_1]}^\perp
v_2/\|P_{[v_1]}^\perp v_2\|$, etc. Then if $H_j := [u_j]$, computer
calculation shows the impossibility of complete coupling.

Now let $G$ be a Cayley graph of a finitely generated infinite group $\Gamma$.
Analogies with percolation (see, e.g., \ref b.Lyons:phase/ for a review) and
with minimal spanning forests (see \ref b.LPS:msf/)
suggest the following possibilities.
Let $H$ be a $\Gamma$-invariant subspace of $\ell^2(\edges)$.

\beginbullets

If $H$ is a proper subspace of the star space $\STAR$ of
$G$, then all connected components of $\ba$ are finite $\P^H$-a.s.

If $\STAR \subsetneq H \subsetneq \CYCLE^\perp$, where $\CYCLE$ is the cycle
space of $G$, then $\ba$ has infinitely many (infinite) components $\P^H$-a.s.

If $\CYCLE^\perp \subsetneq H$, then $\ba$ has a single (infinite) component
$\P^H$-a.s. 

\endbullets

The second statement is shown to be true by \ref b.Lyons:betti/.
% Pf: $\E^H[\deg_\ba o] = \dim_\Gamma H$. By ergodicity, the number of
% components is constant. If finite, say, $m$, choose $H = H_0 \subsetneq H_1
% \subsetneq H_2 \subsetneq \cdots \subsetneq H_{m} = \CYCLE^\perp$, all
% invariant. This can be done using vN dimension.
% By domination, the number of components is strictly smaller with
% respect to $\P^{H_{i+1}}$ that it is with respect to $\P^{H_i}$. This is a
% contradiction.
We do not know whether the others are true.
However, when $H \subsetneq \STAR$, the expected
degree of a vertex with respect to $\P^H$ is less than 2, its expected degree
in the $\wsf$ (\BLPSusf).
By Theorem 6.1 of \BLPSgip, it follows that $\ba$ has infinitely many finite
components $\P^H$-a.s.
It is shown in \ref b.Lyons:betti/ that the last bulleted statement above 
implies a positive answer to an important question of
\ref b.Gaboriau:invariants/, showing that the cost of $\Gamma$ is equal to 1
plus the first $\ell^2$-Betti number of $\Gamma$.
\comment{Mention not vacuous questions? That comes from vN dim.}

\comment{
Are minimal base measures interesting (analogous to minimal spanning
forests)? When they are ordered, they have a natural probability measure on
them which is the stationary measure for an analogue of the Tsetlin
library; see Brown (2000).
}

%There should be some theorem of the following flavor.  Let $v_1,..,v_k$ be
%boundary vertices.  Then they are all in different components a.s.\ iff the
%image of $H$ under the map $\theta \mapsto (\theta . star(v_1), ... \theta .
%star(v)_k)$ has dimension $n-k$, where $n$ is the total number of boundary
%vertices, I think.

\bsection{Open Questions: Surface Graphs and CW-Complexes}{s.applic:CW}

%One can consider the process corresponding to the subspace $\nHD$ itself;
%this will be the ``difference'' between the free and the wired spanning
%forests if our coupling is successful. At the moment, this process is
%entirely mysterious to us; yet it is certainly a ``natural'' process in
%that it is invariant under automorphisms.

For graphs $G = (\vertex, \edge)$, 
one need not restrict oneself to subspaces $H$ of $\ell^2(\edge)$
that give spanning
trees or forests. For example, if $G$ is a graph that is embedded on a
surface (such as a punctured plane or the 2-torus),
let $H$ be the orthocomplement of the boundaries (i.e., of the
image of the boundary operator $\bd_2$).
In this case, the measure $\P^H$ is the uniform measure on 
maximal subgraphs that do not contain any boundary
in the sense that no linear
combination of the edges is a boundary.
(This alternative description is proved by the considerations of the
following paragraph.)
Properties of $\P^H$ are worth investigation.
A particular example arises as follows.
If a graph is embedded on a torus and one takes a uniform
spanning tree of the graph, then its complement on the dual graph relative to
the surface contains only noncontractible cycles. The distribution on
homology is worth investigation.
Can we calculate the distribution of the (unsigned) homology basis it gives?
Does
it have a limit as the mesh of the graph tends to 0?  Presumably the limit
does exist and is conformally invariant, but it would greatly help if one
could find a way to generate the homology basis directly without generating
the entire subgraph, such as via some algorithm analogous to those of
Aldous/Broder or Wilson.
%For these measures and others, perhaps there is an analogue of
%Wilson's algorithm, either via random walk or via ``cycle-popping''.

For another class of examples, consider a finite CW-complex $K$ of
dimension $d$.
Given $0 < k \le d$, the representation of the matroid corresponding to the
matrix (with respect to the usual cellular bases) of the boundary operator
$\bd_{k}$ from $k$-chains to $(k-1)$-chains yields a probability measure
$\pd{k}$ on the set of maximal $k$-subcomplexes $L$ of $K$ with $H_{k} (L;\,
\Q) = 0$. 
For example, $\pd1$ is the uniform spanning tree on the 1-skeleton of $K$.
In general, as shown by \ref b.Lyons:betti/, the probability of such a
subcomplex $L$ is proportional to the square of the order of the torsion
subgroup of $H_{k-1}(L;\, \Z)$, the $(k-1)$-dimensional homology group of the
subcomplex $L$.
When $K$ is a simplex, \ref b.Kalai/ showed that the number of 
maximal $\Q$-acyclic $k$-subcomplexes of $K$ counted with these weights is
$$
n^{\raise5pt\hbox{${\displaystyle n-2 \choose \displaystyle k}$}}
\,,
$$
thereby generalizing Cayley's theorem.
In case $K$ is infinite and locally finite, one can take free and wired
limits analogous to the $\fusf$ and the $\wusf$.
In \ref b.Lyons:betti/, it is shown that in any amenable transitive
contractible complex, these free and wired measures agree, which provides a
new proof\comment{and extension} of a theorem of \ref b.CheegerGromov/ 
concerning the vanishing of $\ell^2$-Betti numbers.

Now specialize to the natural $d$-dimensional CW-complex determined by
the hyperplanes of $\R^d$ passing through points of $\Z^d$ and parallel to
the coordinate hyperplanes (so the 0-cells are the points of
$\Z^d$). 
In \ref b.Lyons:betti/, it is shown that the
$\pd{k}$-probability that a given $k$-cell belongs to the random
$k$-subcomplex is $k/d$.
But many questions are open.
Among the most important are the following two:

\beginbullets

What is the $(k-1)$-dimensional (co)homology of the
$k$-subcomplex? In the case $k=1$ of spanning forests, this asks how many
trees there are, the question answered by
\ref b.Pemantle:ust/.

If one takes the 1-point compactification of the subcomplex, what is the
$k$-dimensional (co)homology? In the case of spanning forests, this asks
how many ends there are in the tree(s), the 
question answered partially by
\ref b.Pemantle:ust/ and completely by \BLPSusf.

\endbullets

Note that by translation-invariance of (co)homology and ergodicity of $\pd k$,
we have that the values of the (co)homology groups are constants a.s.

The two questions above are interesting even for rational (co)homology. It
then follows trivially from the Alexander duality theorem and the results
of \ref b.Pemantle:ust/ and \BLPSusf\ 
that for $k=d-1$, we have $H_{k-1}(L;\, \Q) = 0$ $\pd k$-a.s., while
$\pd k$-a.s.\ $H^k(L \cup \infty;\, \Q)$ is $0$ for $2 \le d \le 4$ and is
(naturally isomorphic to)
an infinite direct product of $\Q$ for $d \ge 5$ (so the homology is the
infinite direct sum of $\Q$).
It also follows from the Alexander duality theorem and from equality of free
and wired limits that if $d=2k$, then the a.s.\ values of
$H^k(L \cup \infty;\, \Q)$ and $H_{k-1}(L;\, \Q)$ are the same (naturally
isomorphic), so that the two bulleted questions above are dual in that case.
Since even the finite complexes can have nontrivial integral homology, it is
probably more interesting to examine the quotient of integral cohomology by
integral cohomology with compact support.
%, or else simply rational (co)homology.
In the present case, the finite complexes have finite groups $H_{k-1}(L;\,
\Z)$, so we might simply ask about the Betti numbers in the infinite limit.

Many matroids, of course, are not representable.
For them, the above theory gives no measure on the bases.
A first clue of how to define an interesting measure nevertheless
comes from the following observation.
Suppose we choose a uniform spanning tree from a graph that has $n$ vertices.
If we then choose an edge uniformly from the tree, the chance of picking $e$
is $Y(e, e)/(n-1)$ by Kirchhoff's Theorem, where $Y$ is the transfer current
matrix.
Therefore $\sum_{e \in E} Y(e, e) = n-1$, a theorem of \ref b.Foster/
on electrical networks.
One might thus expect something interesting from the measures we now
introduce, even when specialized to a graphical matroid.

A second clue is that every matroid $\M = (E, \B)$ has naturally associated to
it a simplicial complex $K_\M$ 
%the {\bf matroid complex}
formed from its independent sets, where a subset
of $E$ is called {\bf independent} if it lies in some base.

We now see how to define a natural probability measure on $\B$.
Namely, let $r$ be the
rank of $\M$, which is one more than the dimension of $K_\M$.
The boundary operator $\bd_{r-1}$ for $K_\M$ gives, as above, a
probability measure $\pd{r-1}$ on collections of bases, and then one may choose
uniformly an element of such a collection to obtain, finally, a probability
measure $\pcx{r-1}$ on $\B$.
This is not just a complicated way of defining the uniform measure on $\B$,
yet it is a measure that is invariant under automorphisms of the matroid.
Because of this, it gives a new measure even for the graphic matroid, i.e.,
a new measure on spanning trees of a graph that, although not uniform,
is invariant under automorphisms of the graph, as well as under all
matroid automorphisms.
Thus, this new measure reflects more of the structure of the graph and of
the matroid than does the uniform measure.
If one does not use the top dimension $r-1$ but a dimension $k < r-1$, then one 
obtains an automorphism-invariant probability measure $\pcx{k}$
on the independent sets of cardinality $k$.

\beginbullets

How are these measures $\pcx{k}$ related to each other?

How do they behave under the standard matroid operations?

Can one describe the measures more explicitly and directly for graphic
matroids?

\endbullets

\bsection{Open Questions: Dynamical Systems}{s.applic:dyn}

Let $\Gamma$ be a countable infinite discrete abelian group, such as $\Z^n$.
Let $\widehat \Gamma$ be the group dual to $\Gamma$, a compact group
equal to $\R^n/\Z^n$ when $\Gamma = \Z^n$.
If $f : \widehat \Gamma \to [0, 1]$ is a measurable function,
then multiplication by $f$ is a
positive contraction on $L^2_\C(\widehat\Gamma, \lambda)$, where $\lambda$
is unit Haar measure.
Since $L^2_\C(\widehat\Gamma, \lambda)$ is isomorphic to $\ell^2(\Gamma; \C)$
via the Fourier transform, there is an associated probability measure $\P^f$
on $2^{\Gamma}$.
As an easy example, if $f$ is a constant, $p$, then $\P^f$ is just the
Bernoulli($p$) process on $\Gamma$.
In general, the measure $\P^f$ is invariant under the natural $\Gamma$
action and has a trivial full tail $\sigma$-field by \ref t.tail/.
\ref b.LS:dyn/ have shown that if $\Gamma = \Z^n$, then 
for any $f$, we also have
that the dynamical system $(2^\Gamma, \P^f, \Gamma)$ is a
Bernoulli shift, i.e., is isomorphic to an i.i.d.\ process.
Therefore, by Ornstein's theorem (and its generalizations, see \ref
b.KatzWeiss/, \ref b.Conze/, \ref b.Thouvenot/, and \ref b.OrnWeiss/),
it is characterized up to isomorphism by its entropy.
%What about amenable $\Gamma$?

\beginbullets

What is the entropy of the dynamical system $(2^\Gamma, \P^f, \Gamma)$?

We conjecture that entropy is concave.
In other words, if $h(f)$ denotes the entropy of the
dynamical system $(2^\Gamma, \P^f, \Gamma)$, then
$h\big({(f+g)/2}\big) \ge \big(h(f) + h(g)\big)/2$ for all $f$ and $g$.
This is a corollary of \ref g.concave/, and even from a
restricted version of \ref e.concave/ that assumes that both $Q_1$ and
$Q_2$ are Toeplitz matrices.
% circulant?? This is by looking at an $n$-step step function approx to $f$
% and F.T. on $\Z_n$.
%$Q_1 Q_2 = Q_2 Q_1$ and $Q_1$, $Q_2$ are Toeplitz matrices.

If $0 \le f \le g \le 1$, then 
there is a monotone coupling of $\P^f$ and $\P^g$ by \ref t.Q/. Can an explicit
monotone coupling be given? For example, if $f = pg$, where $p$ is
a constant, then this
can be done by using the fact that $\P^{pg}$ has the same law as the
pointwise minimum of independent processes $\P^p$ and $\P^g$.

Consider the case $\Gamma = \Z$.
Note that translation and flip of $f$ yields the same measure $\P^{f}$,
even though $f$ changes. Does $\P^{f}$ determine $f$ up to translation and
flip? 

\endbullets

Additional questions concerning these systems appear in \ref b.LS:dyn/.

\comment{
In the special case where 
$f$ is the indicator of a set $A \subseteq \widehat\Gamma$,
we obtain an action of $\Gamma$ on
the subspace $L^2(A)$ of functions that are 0 outside $A$.
One may regard this action of $\Gamma$ on $L^2(A)$ as a representation of
$\Gamma$. This suggests considering the dynamical system on $2^\Gamma$
associated to any subrepresentation of the regular representation of a
countable group $\Gamma$ by the probability measure $\P^H$ given by the
invariant subspace $H \subseteq \ell^2(\Gamma)$ that determines the
subrepresentation. How are properties of the subrepresentation reflected in
properties of the dynamical system? For example, one checks easily that
the von Neumann dimension of the subrepresentation is equal to the
one-dimensional marginal of $\P^H$, i.e., the probability that any given
element of $\Gamma$ belongs to $\ba$.
questions.}

\medbreak
\noindent {\bf Acknowledgements.}\enspace
I am indebted to Oded Schramm for many conversations at the beginning of
this project in Fall 1997; in particular, it was at this time that we
developed the idea of using exterior algebra as a convenient representation
of determinantal probability measures. His influence and contributions
appear in several
other places as well, and he made a number of useful comments on a near-final
draft.
I thank Hari Bercovici and L\'aszl\'o Lov\'asz for some helpful discussions.
I am grateful to Alexander Soshnikov for several references and to Alexander
Ulanovskii for useful remarks.

%\bibfile{\jobname}
\def\noop#1{\relax}
\input \jobname.bbl

\filbreak
\begingroup
\eightpoint\sc
\parindent=0pt\baselineskip=10pt

Department of Mathematics,
Indiana University,
Bloomington, IN 47405-5701, USA
\emailwww{rdlyons@indiana.edu}
{http://mypage.iu.edu/\string~rdlyons/}

and

School of Mathematics,
Georgia Institute of Technology,
Atlanta, GA 30332-0160

\endgroup

\tracingstats=1


\def\cprime{$'$} \def\cprime{$'$}
\def\temp{\let\linkit=\linkyear \apaliketrue}
\temp
\ifcitationgeneration\immediate\write\labelfile{\sanitize\temp}\fi
\def\startreferences{
 \vskip0pt plus.3\vsize \penalty -150 \vskip0pt
 plus-.3\vsize \bigskip\bigskip \vskip \parskip
 \begingroup\baselineskip=12pt\frenchspacing
 \bibliographytitle
 \vskip12pt\parindent=0pt
 \def\and{{\rm and}}
 \def\em{\it}
 \def\newblock{\hskip .11em plus.33em minus.07em}
 \def\bibauthor##1{{\sc ##1}}
 \def\bibitem[##1]##2
 {\htmlanchor{##2}{}\RefLabel{##2}[##1]\hangindent=.8cm\hangafter=1}
 }
\def\endreferences{\bigskip\bigskip\endgroup}
\ifundefined{bibstylemodification}\relax\else\bibstylemodification\fi
\startreferences

\bibitem[Aldous (1990)]{MR91h:60013}
\bibauthor{Aldous, D.J.} (1990).
\newblock The random walk construction of uniform spanning trees and uniform
  labelled trees.
\newblock {\em SIAM J. Discrete Math.} {\bf 3}, 450--465.

\bibitem[Alon and Spencer (2001)]{AS:book}
\bibauthor{Alon, N. \and{} Spencer, J.H.} (2001).
\newblock {\em The Probabilistic Method}.
\newblock John Wiley \& Sons Inc., New York, second edition.

\bibitem[Benjamini, Lyons, Peres, and Schramm (1999)]{MR1675890}
\bibauthor{Benjamini, I., Lyons, R., Peres, Y., \and{} Schramm, O.} (1999).
\newblock Group-invariant percolation on graphs.
\newblock {\em Geom. Funct. Anal.} {\bf 9}, 29--66.

\bibitem[Benjamini, Lyons, Peres, and Schramm (2001)]{BLPSusf}
\bibauthor{Benjamini, I., Lyons, R., Peres, Y., \and{} Schramm, O.} (2001).
\newblock Uniform spanning forests.
\newblock {\em Ann. Probab.} {\bf 29}, 1--65.

\bibitem[van~den Berg and Kesten (1985)]{MR87b:60027}
\bibauthor{van~den Berg, J. \and{} Kesten, H.} (1985).
\newblock Inequalities with applications to percolation and reliability.
\newblock {\em J. Appl. Probab.} {\bf 22}, 556--569.

\bibitem[Beurling and Malliavin (1967)]{MR35:654}
\bibauthor{Beurling, A. \and{} Malliavin, P.} (1967).
\newblock On the closure of characters and the zeros of entire functions.
\newblock {\em Acta Math.} {\bf 118}, 79--93.

\bibitem[Borodin (2000)]{MR2001h:20017}
\bibauthor{Borodin, A.} (2000).
\newblock Characters of symmetric groups, and correlation functions of point
  processes.
\newblock {\em Funktsional. Anal. i Prilozhen.} {\bf 34}, 12--28, 96.
\newblock English translation: {\it Funct. Anal. Appl.} {\bf 34} (2000), no. 1,
  10--23.

\bibitem[Borodin, Okounkov, and Olshanski (2000)]{MR2001g:05103}
\bibauthor{Borodin, A., Okounkov, A., \and{} Olshanski, G.} (2000).
\newblock Asymptotics of {P}lancherel measures for symmetric groups.
\newblock {\em J. Amer. Math. Soc.} {\bf 13}, 481--515 (electronic).

\bibitem[Borodin and Olshanski (2000)]{MR2001k:33031}
\bibauthor{Borodin, A. \and{} Olshanski, G.} (2000).
\newblock Distributions on partitions, point processes, and the hypergeometric
  kernel.
\newblock {\em Comm. Math. Phys.} {\bf 211}, 335--358.

\bibitem[Borodin and Olshanski (2001)]{MR2002m:82026}
\bibauthor{Borodin, A. \and{} Olshanski, G.} (2001).
\newblock {$z$}-measures on partitions, {R}obinson-{S}chensted-{K}nuth
  correspondence, and {$\beta=2$} random matrix ensembles.
\newblock In Bleher, P. \and{} Its, A., editors, {\em Random Matrix Models and
  Their Applications}, volume 40 of {\em Math. Sci. Res. Inst. Publ.}, pages
  71--94. Cambridge Univ. Press, Cambridge.

\bibitem[Borodin and Olshanski (2002)]{BO:harmonic}
\bibauthor{Borodin, A. \and{} Olshanski, G.} (2002).
\newblock Harmonic analysis on the infinite-dimensional unitary group and
  determinantal point processes.
\newblock Preprint.

\bibitem[Bourgain and Tzafriri (1987)]{MR89a:46035}
\bibauthor{Bourgain, J. \and{} Tzafriri, L.} (1987).
\newblock Invertibility of ``large'' submatrices with applications to the
  geometry of {B}anach spaces and harmonic analysis.
\newblock {\em Israel J. Math.} {\bf 57}, 137--224.

\bibitem[Broder (1989)]{Broder}
\bibauthor{Broder, A.} (1989).
\newblock Generating random spanning trees.
\newblock In {\em 30th Annual Symposium on Foundations of Computer Science
  (Research Triangle Park, North Carolina)}, pages 442--447, New York. IEEE.

\bibitem[Brooks, Smith, Stone, and Tutte (1940)]{MR2:153d}
\bibauthor{Brooks, R.L., Smith, C.A.B., Stone, A.H., \and{} Tutte, W.T.}
  (1940).
\newblock The dissection of rectangles into squares.
\newblock {\em Duke Math. J.} {\bf 7}, 312--340.

\bibitem[Burton and Pemantle (1993)]{MR94m:60019}
\bibauthor{Burton, R.M. \and{} Pemantle, R.} (1993).
\newblock Local characteristics, entropy and limit theorems for spanning trees
  and domino tilings via transfer-impedances.
\newblock {\em Ann. Probab.} {\bf 21}, 1329--1371.

\bibitem[Cheeger and Gromov (1986)]{MR87i:58161}
\bibauthor{Cheeger, J. \and{} Gromov, M.} (1986).
\newblock ${L}\sb 2$-cohomology and group cohomology.
\newblock {\em Topology} {\bf 25}, 189--215.

\bibitem[Choe, Oxley, Sokal, and Wagner (2003)]{COSW}
\bibauthor{Choe, Y.B., Oxley, J., Sokal, A., \and{} Wagner, D.} (2003).
\newblock Homogeneous multivariate polynomials with the half-plane property.
\newblock {\em Adv. in Appl. Math.}
\newblock To appear.

\bibitem[Conrey (2003)]{Conrey}
\bibauthor{Conrey, J.B.} (2003).
\newblock The {R}iemann hypothesis.
\newblock {\em Notices Amer. Math. Soc.} {\bf 50}, 341--353.

\bibitem[Conway (1990)]{MR91e:46001}
\bibauthor{Conway, J.B.} (1990).
\newblock {\em A Course in Functional Analysis}.
\newblock Springer-Verlag, New York, second edition.

\bibitem[Conze (1972/73)]{MR49:534}
\bibauthor{Conze, J.P.} (1972/73).
\newblock Entropie d'un groupe ab\'elien de transformations.
\newblock {\em Z. Wahrscheinlichkeitstheorie und Verw. Gebiete} {\bf 25},
  11--30.

\bibitem[Daley and Vere-Jones (1988)]{MR90e:60060}
\bibauthor{Daley, D.J. \and{} Vere-Jones, D.} (1988).
\newblock {\em An Introduction to the Theory of Point Processes}.
\newblock Springer-Verlag, New York.

\bibitem[Diaconis (2003)]{MR1962294}
\bibauthor{Diaconis, P.} (2003).
\newblock Patterns in eigenvalues: the 70th {J}osiah {W}illard {G}ibbs lecture.
\newblock {\em Bull. Amer. Math. Soc. (N.S.)} {\bf 40}, 155--178 (electronic).

\bibitem[Dubhashi and Ranjan (1998)]{MR99k:60048}
\bibauthor{Dubhashi, D. \and{} Ranjan, D.} (1998).
\newblock Balls and bins: a study in negative dependence.
\newblock {\em Random Structures Algorithms} {\bf 13}, 99--124.

\bibitem[Dyson (1962)]{MR26:1113}
\bibauthor{Dyson, F.J.} (1962).
\newblock Statistical theory of the energy levels of complex systems. {III}.
\newblock {\em J. Mathematical Phys.} {\bf 3}, 166--175.

\bibitem[Feder and Mihail (1992)]{FedMih}
\bibauthor{Feder, T. \and{} Mihail, M.} (1992).
\newblock Balanced matroids.
\newblock In {\em Proceedings of the {T}wenty-{F}ourth {A}nnual {A}{C}{M}
  {S}ymposium on {T}heory of {C}omputing}, pages 26--38, New York. Association
  for Computing Machinery (ACM).
\newblock Held in Victoria, BC, Canada.

\bibitem[Foster (1948)]{MR10:662a}
\bibauthor{Foster, R.M.} (1948).
\newblock The average impedance of an electrical network.
\newblock In {\em Reissner Anniversary Volume, Contributions to Applied
  Mechanics}, pages 333--340. J. W. Edwards, Ann Arbor, Michigan.
\newblock Edited by the Staff of the Department of Aeronautical Engineering and
  Applied Mechanics of the Polytechnic Institute of Brooklyn.

\bibitem[Fulton and Harris (1991)]{MR93a:20069}
\bibauthor{Fulton, W. \and{} Harris, J.} (1991).
\newblock {\em Representation Theory: A First Course}.
\newblock Springer-Verlag, New York.
\newblock Readings in Mathematics.

\bibitem[Gaboriau (2002)]{MR1953191}
\bibauthor{Gaboriau, D.} (2002).
\newblock Invariants {$l\sp 2$} de relations d'\'equivalence et de groupes.
\newblock {\em Publ. Math. Inst. Hautes \'Etudes Sci.}, 93--150.

\bibitem[Georgii (1988)]{MR89k:82010}
\bibauthor{Georgii, H.O.} (1988).
\newblock {\em Gibbs Measures and Phase Transitions}.
\newblock Walter de Gruyter \& Co., Berlin-New York.

\bibitem[H{\"a}ggstr{\"o}m (1995)]{MR97b:60170}
\bibauthor{H{\"a}ggstr{\"o}m, O.} (1995).
\newblock Random-cluster measures and uniform spanning trees.
\newblock {\em Stochastic Process. Appl.} {\bf 59}, 267--275.

\bibitem[Halmos (1982)]{MR84e:47001}
\bibauthor{Halmos, P.R.} (1982).
\newblock {\em A {H}ilbert {S}pace {P}roblem {B}ook}.
\newblock Springer-Verlag, New York, second edition.
\newblock Encyclopedia of Mathematics and its Applications, 17.

\bibitem[Heicklen and Lyons (2003)]{HeicklenLyons}
\bibauthor{Heicklen, D. \and{} Lyons, R.} (2003).
\newblock Change intolerance in spanning forests.
\newblock {\em J. Theoret. Probab.} {\bf 16}, 47--58.

\bibitem[Johansson (2001)]{MR1826414}
\bibauthor{Johansson, K.} (2001).
\newblock Discrete orthogonal polynomial ensembles and the {P}lancherel
  measure.
\newblock {\em Ann. of Math. (2)} {\bf 153}, 259--296.

\bibitem[Johansson (2002)]{MR1900323}
\bibauthor{Johansson, K.} (2002).
\newblock Non-intersecting paths, random tilings and random matrices.
\newblock {\em Probab. Theory Related Fields} {\bf 123}, 225--280.

\bibitem[Kalai (1983)]{MR85a:55006}
\bibauthor{Kalai, G.} (1983).
\newblock Enumeration of ${\Q}$-acyclic simplicial complexes.
\newblock {\em Israel J. Math.} {\bf 45}, 337--351.

\bibitem[Katznelson and Weiss (1972)]{MR47:5227}
\bibauthor{Katznelson, Y. \and{} Weiss, B.} (1972).
\newblock Commuting measure-preserving transformations.
\newblock {\em Israel J. Math.} {\bf 12}, 161--173.

\bibitem[Kirchhoff (1847)]{Kirch}
\bibauthor{Kirchhoff, G.} (1847).
\newblock Ueber die {A}ufl\"osung der {G}leichungen, auf welche man bei der
  {U}ntersuchung der linearen {V}ertheilung galvanischer {S}tr\"ome gef\"uhrt
  wird.
\newblock {\em Ann. Phys. und Chem.} {\bf 72}, 497--508.

\bibitem[Lyons (1998)]{MR1630412}
\bibauthor{Lyons, R.} (1998).
\newblock A bird's-eye view of uniform spanning trees and forests.
\newblock In Aldous, D. \and{} Propp, J., editors, {\em Microsurveys in
  Discrete Probability}, volume 41 of {\em DIMACS Series in Discrete
  Mathematics and Theoretical Computer Science}, pages 135--162. Amer. Math.
  Soc., Providence, RI.
\newblock Papers from the workshop held as part of the Dimacs Special Year on
  Discrete Probability in Princeton, NJ, June 2--6, 1997.

\bibitem[Lyons (2000)]{MR2001c:82028}
\bibauthor{Lyons, R.} (2000).
\newblock Phase transitions on nonamenable graphs.
\newblock {\em J. Math. Phys.} {\bf 41}, 1099--1126.
\newblock Probabilistic techniques in equilibrium and nonequilibrium
  statistical physics.

\bibitem[Lyons (2003)]{Lyons:betti}
\bibauthor{Lyons, R.} (2003).
\newblock Random complexes and $\ell^2$-{B}etti numbers.
\newblock In preparation.

\bibitem[Lyons, Peres, and Schramm (2003)]{LPS:msf}
\bibauthor{Lyons, R., Peres, Y., \and{} Schramm, O.} (2003).
\newblock Minimal spanning forests.
\newblock In preparation.

\bibitem[Lyons and Steif (2003)]{LyonsSteif:dyn}
\bibauthor{Lyons, R. \and{} Steif, J.E.} (2003).
\newblock Stationary determinantal processes: Phase multiplicity,
  {B}ernoullicity, entropy, and domination.
\newblock {\em Duke Math. J.}
\newblock To appear.

\bibitem[Macchi (1975)]{MR52:1876}
\bibauthor{Macchi, O.} (1975).
\newblock The coincidence approach to stochastic point processes.
\newblock {\em Advances in Appl. Probability} {\bf 7}, 83--122.

\bibitem[Maurer (1976)]{MR52:13452}
\bibauthor{Maurer, S.B.} (1976).
\newblock Matrix generalizations of some theorems on trees, cycles and cocycles
  in graphs.
\newblock {\em SIAM J. Appl. Math.} {\bf 30}, 143--148.

\bibitem[Mehta (1991)]{MR92f:82002}
\bibauthor{Mehta, M.L.} (1991).
\newblock {\em Random Matrices}.
\newblock Academic Press Inc., Boston, MA, second edition.

\bibitem[Morris (2002)]{Morris}
\bibauthor{Morris, B.} (2002).
\newblock The components of the wired spanning forest are recurrent.
\newblock {\em Probab. Theory Related Fields}.
\newblock To appear.

\bibitem[Newman (1984)]{MR86i:60072}
\bibauthor{Newman, C.M.} (1984).
\newblock Asymptotic independence and limit theorems for positively and
  negatively dependent random variables.
\newblock In Tong, Y.L., editor, {\em Inequalities in Statistics and
  Probability}, pages 127--140. Inst. Math. Statist., Hayward, CA.
\newblock Proceedings of the symposium held at the University of Nebraska,
  Lincoln, Neb., October 27--30, 1982.

\bibitem[Okounkov (2001)]{MR2002f:60019}
\bibauthor{Okounkov, A.} (2001).
\newblock Infinite wedge and random partitions.
\newblock {\em Selecta Math. (N.S.)} {\bf 7}, 57--81.

\bibitem[Okounkov and Reshetikhin (2001)]{OkoResh:schur}
\bibauthor{Okounkov, A. \and{} Reshetikhin, N.} (2001).
\newblock Correlation function of {S}chur process with application to local
  geometry of a random 3-dimensional {Y}oung diagram.
\newblock Preprint.

\bibitem[Ornstein and Weiss (1987)]{MR88j:28014}
\bibauthor{Ornstein, D.S. \and{} Weiss, B.} (1987).
\newblock Entropy and isomorphism theorems for actions of amenable groups.
\newblock {\em J. Analyse Math.} {\bf 48}, 1--141.

\bibitem[Oxley (1992)]{MR94d:05033}
\bibauthor{Oxley, J.G.} (1992).
\newblock {\em Matroid Theory}.
\newblock The Clarendon Press Oxford University Press, New York.

\bibitem[Pemantle (1991)]{MR92g:60014}
\bibauthor{Pemantle, R.} (1991).
\newblock Choosing a spanning tree for the integer lattice uniformly.
\newblock {\em Ann. Probab.} {\bf 19}, 1559--1574.

\bibitem[Pemantle (2000)]{MR2001g:62039}
\bibauthor{Pemantle, R.} (2000).
\newblock Towards a theory of negative dependence.
\newblock {\em J. Math. Phys.} {\bf 41}, 1371--1390.
\newblock Probabilistic techniques in equilibrium and nonequilibrium
  statistical physics.

\bibitem[Propp and Wilson (1998)]{MR99g:60116}
\bibauthor{Propp, J.G. \and{} Wilson, D.B.} (1998).
\newblock How to get a perfectly random sample from a generic {M}arkov chain
  and generate a random spanning tree of a directed graph.
\newblock {\em J. Algorithms} {\bf 27}, 170--217.
\newblock 7th Annual ACM-SIAM Symposium on Discrete Algorithms (Atlanta, GA,
  1996).

\bibitem[Redheffer (1972)]{MR48:801}
\bibauthor{Redheffer, R.} (1972).
\newblock Two consequences of the {B}eurling-{M}alliavin theory.
\newblock {\em Proc. Amer. Math. Soc.} {\bf 36}, 116--122.

\bibitem[Redheffer (1977)]{MR56:5852}
\bibauthor{Redheffer, R.M.} (1977).
\newblock Completeness of sets of complex exponentials.
\newblock {\em Advances in Math.} {\bf 24}, 1--62.

\bibitem[Seip and Ulanovskii (1997)]{MR97g:42003}
\bibauthor{Seip, K. \and{} Ulanovskii, A.M.} (1997).
\newblock The {B}eurling-{M}alliavin density of a random sequence.
\newblock {\em Proc. Amer. Math. Soc.} {\bf 125}, 1745--1749.

\bibitem[Shao (2000)]{MR2001g:60077}
\bibauthor{Shao, Q.M.} (2000).
\newblock A comparison theorem on moment inequalities between negatively
  associated and independent random variables.
\newblock {\em J. Theoret. Probab.} {\bf 13}, 343--356.

\bibitem[Shao and Su (1999)]{MR2000g:60052}
\bibauthor{Shao, Q.M. \and{} Su, C.} (1999).
\newblock The law of the iterated logarithm for negatively associated random
  variables.
\newblock {\em Stochastic Process. Appl.} {\bf 83}, 139--148.

\bibitem[Shirai and Takahashi (2000)]{ShiTak:ann}
\bibauthor{Shirai, T. \and{} Takahashi, Y.} (2000).
\newblock Fermion process and {F}redholm determinant.
\newblock In H.G.W.~Begehr, R.G. \and{} Kajiwara, J., editors, {\em Proceedings
  of the Second ISAAC Congress}, volume 1, pages 15--23. Kluwer Academic Publ.
\newblock International Society for Analysis, Applications and Computation
  Volume 7.

\bibitem[Shirai and Takahashi (2002)]{ShiTak:I}
\bibauthor{Shirai, T. \and{} Takahashi, Y.} (2002).
\newblock Random point fields associated with certain {F}redholm determinants
  {I}: fermion, {P}oisson and boson point processes.
\newblock Preprint.

\bibitem[Shirai and Takahashi (2003)]{ShiTak:II}
\bibauthor{Shirai, T. \and{} Takahashi, Y.} (2003).
\newblock Random point fields associated with certain {F}redholm determinants
  {II}: fermion shifts and their ergodic properties.
\newblock {\em Ann. Probab.}
\newblock To appear.

\bibitem[Shirai and Yoo (2002)]{MR1942400}
\bibauthor{Shirai, T. \and{} Yoo, H.J.} (2002).
\newblock Glauber dynamics for fermion point processes.
\newblock {\em Nagoya Math. J.} {\bf 168}, 139--166.

\bibitem[Soshnikov (2000a)]{MR1799012}
\bibauthor{Soshnikov, A.} (2000a).
\newblock Determinantal random point fields.
\newblock {\em Uspekhi Mat. Nauk} {\bf 55}, 107--160.

\bibitem[Soshnikov (2000b)]{MR2001m:82006}
\bibauthor{Soshnikov, A.B.} (2000b).
\newblock Gaussian fluctuation for the number of particles in {A}iry, {B}essel,
  sine, and other determinantal random point fields.
\newblock {\em J. Statist. Phys.} {\bf 100}, 491--522.

\bibitem[Strassen (1965)]{MR31:1693}
\bibauthor{Strassen, V.} (1965).
\newblock The existence of probability measures with given marginals.
\newblock {\em Ann. Math. Statist} {\bf 36}, 423--439.

\bibitem[Thomassen (1990)]{MR91d:94029}
\bibauthor{Thomassen, C.} (1990).
\newblock Resistances and currents in infinite electrical networks.
\newblock {\em J. Combin. Theory Ser. B} {\bf 49}, 87--102.

\bibitem[Thouvenot (1972)]{MR47:10145}
\bibauthor{Thouvenot, J.P.} (1972).
\newblock Convergence en moyenne de l'information pour l'action de {${\bf
  Z}\sp{2}$}.
\newblock {\em Z. Wahrscheinlichkeitstheorie und Verw. Gebiete} {\bf 24},
  135--137.

\bibitem[Vershik and Kerov (1981)]{MR84a:22016}
\bibauthor{Vershik, A.M. \and{} Kerov, S.V.} (1981).
\newblock Asymptotic theory of the characters of a symmetric group.
\newblock {\em Funktsional. Anal. i Prilozhen.} {\bf 15}, 15--27, 96.
\newblock English translation: {\it Functional Anal. Appl.} {\bf 15} (1981),
  no. 4, 246--255 (1982).

\bibitem[Welsh (1976)]{MR55:148}
\bibauthor{Welsh, D.J.A.} (1976).
\newblock {\em Matroid Theory}.
\newblock Academic Press [Harcourt Brace Jovanovich Publishers], London.
\newblock L. M. S. Monographs, No. 8.

\bibitem[White (1987)]{MR88g:05048}
\bibauthor{White, N., editor} (1987).
\newblock {\em Combinatorial Geometries}.
\newblock Cambridge University Press, Cambridge.

\bibitem[Whitney (1935)]{Whitney:mat}
\bibauthor{Whitney, H.} (1935).
\newblock On the abstract properties of linear dependence.
\newblock {\em Amer. J. Math.} {\bf 57}, 509--533.

\bibitem[Whitney (1957)]{MR19:309c}
\bibauthor{Whitney, H.} (1957).
\newblock {\em Geometric Integration Theory}.
\newblock Princeton University Press, Princeton, N. J.

\bibitem[Wilson (1996)]{MR1427525}
\bibauthor{Wilson, D.B.} (1996).
\newblock Generating random spanning trees more quickly than the cover time.
\newblock In {\em Proceedings of the {T}wenty-eighth {A}nnual {ACM} {S}ymposium
  on the {T}heory of {C}omputing}, pages 296--303, New York. ACM.
\newblock Held in Philadelphia, PA, May 22--24, 1996.

\bibitem[Zhang (2001)]{MR1854030}
\bibauthor{Zhang, L.X.} (2001).
\newblock Strassen's law of the iterated logarithm for negatively associated
  random vectors.
\newblock {\em Stochastic Process. Appl.} {\bf 95}, 311--328.

\bibitem[Zhang and Wen (2001)]{MR1841627}
\bibauthor{Zhang, L.X. \and{} Wen, J.} (2001).
\newblock A weak convergence for negatively associated fields.
\newblock {\em Statist. Probab. Lett.} {\bf 53}, 259--267.

\endreferences
\bye